\documentclass[11pt]{amsart}
\usepackage[T1]{fontenc}
\usepackage[utf8]{inputenc}
\usepackage{lmodern}
\usepackage{amsmath,amssymb,mathrsfs,mathtools}
\usepackage[a4paper,margin=28mm]{geometry}
\usepackage{microtype}
\usepackage[shortlabels]{enumitem}
\usepackage{booktabs}
\usepackage{tikz}
\usepackage[hidelinks]{hyperref}

\newtheorem{theorem}{Theorem}[section]
\newtheorem{proposition}[theorem]{Proposition}
\newtheorem{lemma}[theorem]{Lemma}
\newtheorem{corollary}[theorem]{Corollary}
\newtheorem{introthm}{Theorem}

\theoremstyle{definition}
\newtheorem{definition}[theorem]{Definition}
\newtheorem{example}[theorem]{Example}
\newtheorem{question}[theorem]{Question}
\theoremstyle{remark}
\newtheorem{remark}[theorem]{Remark}
\numberwithin{equation}{section}

\newcommand{\Z}{\mathbb Z}
\newcommand{\R}{\mathbb R}
\newcommand{\Q}{\mathbb Q}
\newcommand{\C}{\mathbb C}

\newcommand{\ZN}{\Z_N}
\newcommand{\ones}{\boldsymbol{1}}
\newcommand{\Dg}{\operatorname{Dg}}
\newcommand{\Ed}{\operatorname{Ed}}
\DeclareMathOperator{\diag}{diag}
\DeclareMathOperator{\ord}{ord}

\title[Cube tilings and rational trees]{Cube tilings, partition counters,\protect\\
and the Calkin--Wilf and Kepler trees}
\author[J. Grytczuk]{Jaros\l aw Grytczuk}
\address{Faculty of Mathematics and Information Science,
Warsaw University of Technology, 00-662 Warsaw, Poland}
\email{jaroslaw.grytczuk@pw.edu.pl}
\author[A. P. Kisielewicz]{Andrzej P. Kisielewicz}
\address{Instytut Matematyki, Informatyki i Ekonometrii,
Uniwersytet Zielonog\'orski, ul. Podg\'orna 50,
65-246 Zielona G\'ora, Poland}
\email{A.Kisielewicz@wmie.uz.zgora.pl}
\author[K. Przes\l awski]{Krzysztof Przes\l awski}
\address{Instytut Matematyki, Informatyki i Ekonometrii,
Uniwersytet Zielonog\'orski, ul. Podg\'orna 50,
65-246 Zielona G\'ora, Poland}
\email{K.Przeslawski@wmie.uz.zgora.pl}
\date{}
\subjclass[2020]{Primary 52C22, 05B45; Secondary 11B57, 11B37, 20K01}
\keywords{Cube tiling, lattice tiling, cocyclic lattice, Stern's diatomic
sequence, Calkin--Wilf tree, Kepler tree, finite Fourier analysis}

\begin{document}
\begin{abstract}
A tiling of the integer lattice $\Z^d$ by translates of the discrete cube
$K_d=\{0,1\}^d$ cuts every translate $K_d+x$ into subcubes. We study the
function $\ell(x)$ counting these parts. Its sum over every translate of
$K_d$ equals $3^d$, for an arbitrary translation set. For lattice tilings
with a triangular basis we obtain a dimension-reduction formula, and we
describe every cocyclic lattice by a single normalized congruence. Along a
cyclic generator, the counter of a signed bidiagonal lattice is a row of
Stern's diatomic sequence; pairs of counter values realize the
Calkin--Wilf tree, and complementary binary words realize the Kepler tree.
For the binary lattice, the translated counter labelings of the cube carry
generating face partitions, supplied by the tiling itself, which force the
levels of the Kepler tree, and every nontrivial face partition has a unique
lamination direction. Using finite Fourier analysis we prove that the
ordered cyclic profile determines a standard cocyclic cube tiling lattice up
to coordinate reflections. For lattices tiling by discrete boxes we
determine exactly when the analogous statement holds.
\end{abstract}
\maketitle

%======================================================================
\section{Introduction}\label{sec:intro}
%======================================================================
The \emph{Calkin--Wilf tree} \cite{CW} is the infinite binary tree with
root $1/1$ in which each vertex $p/q$ has the two children
\[
 \frac{p}{p+q}\qquad\text{and}\qquad\frac{p+q}{q}.
\]
Every positive rational number occurs in it exactly once, in lowest terms,
and the vertices read level by level are the quotients $s(n)/s(n+1)$ of
consecutive terms of Stern's diatomic sequence \cite{Stern,Northshield}.
It is sequence A002487 in the On-Line Encyclopedia of Integer Sequences
\cite{OEIS}, where many further references can be found.
The \emph{Kepler tree}, in the convention used here, has root $1/2$ and
gives $p/q$ the children $p/(p+q)$ and $q/(p+q)$; it enumerates the
rational numbers strictly between $0$ and $1$, and goes back to Kepler's
\emph{Harmonices Mundi} \cite{Kepler}. Other root conventions and related
constructions of rational trees occur in \cite{ABCM}.

In this paper both trees arise from tilings by cubes. Let
$T\subseteq\Z^d$ be such that the translates $K_d+t$, $t\in T$, of the
discrete cube $K_d=\{0,1\}^d$ partition $\Z^d$. For every $x\in\Z^d$ the
tiles cut the translate $K_d+x$ into subcubes; we call their number
$\ell_T(x)$ the \emph{partition counter}. Geometrically
(Section~\ref{sec:black}), if we place a closed black unit cube at every
point of $T$, then $\ell_T(x)$ is the number of black cubes touching the unit cube
$x+[0,1]^d$; see Figure~\ref{fig:black}. The first observation is the
identity
\[
 \sum_{a\in K_d+z}\ell_T(a)=3^d\qquad(z\in\Z^d),
\]
valid for every tiling (Theorem~\ref{thm:sum}).

Local observations of discrete cube tilings can also reveal probabilistic
phenomena. In a related paper, the second author \cite{Kisielewicz} studies
counts of black unit cubes in translated observation windows of side
length~$3$. The paper also gives three proofs of the identity above in its
equivalent window-counting formulation; the correspondence between the two
counters is explained in Section~\ref{sec:black}.

For a lattice tiling, the counter is periodic, and the Minkowski--Haj\'os
theorem provides a triangular basis. This yields a recursion in the
dimension (Theorem~\ref{thm:recursion}). When the quotient
$\Z^d/\Lambda$ is cyclic, the lattice is the set of solutions of a single
congruence $\sum_iw_ix_i\equiv0\pmod{2^d}$ with uniquely normalized
weights (Theorem~\ref{thm:weights}), and the counter is determined by its
values along the generator $e_d$, the \emph{cyclic profile}.

\begin{introthm}[Theorems~\ref{thm:stern}, \ref{thm:kepler} and
Proposition~\ref{prop:binarycw}]\label{thmA}
Let $\Lambda_d^b=\{x\in\Z^d:\sum_{i=1}^d2^{d-i}x_i\equiv0\pmod{2^d}\}$
and let $\ell$ be its partition counter. Then
$\ell(ke_d)=s(2^d+k)$ for $0\le k<2^d$. For a binary word $v$ put
$\pi(v)=\ell_{\Lambda^b_{|v|}}(v)$, with $\pi(\varnothing)=1$, and let $\bar v$ be
the complementary word. Then $v\mapsto\pi(v)/\pi(\bar v)$ is the
Calkin--Wilf tree, the word $v$ being read as a path from the root. Its
restriction to the words ending in $0$, with $v$ having the children $v0$
and $\bar v0$, is the Kepler tree.
\end{introthm}

The same profile occurs for a whole family of signed bidiagonal lattices,
and for those lattices the Calkin--Wilf tree can be realized directly on
the representative cubes (Theorem~\ref{thm:cw}). The Kepler tree also has a
geometric counterpart. In Section~\ref{sec:keplercubes} we formulate local
conditions on a labeling of the vertices of $K_d$, a sum rule along the
diagonals of a partition of a face, and a matching between those diagonals
and edges, which make no reference to rational trees but force the edge
labels to be the levels of the Kepler tree. All translated counter
labelings $u\mapsto\ell(x+u)$ of the binary lattice satisfy these
conditions, with the required partitions supplied by the tiling itself
(Theorem~\ref{thm:allcubes}). Moreover, every such partition is laminated
in a unique direction, and the numbers of parts on its two sides realize
the Calkin--Wilf levels (Theorem~\ref{thm:lamination}).

The second half of the paper asks how much of a lattice is remembered by
its cyclic profile. Coordinate reflections preserve the profile, and the
converse holds.

\begin{introthm}[Theorem~\ref{thm:main}]\label{thmB}
Two standard cocyclic cube tiling lattices in $\Z^d$ have the same ordered
cyclic profile $(\ell_\Lambda(ke_d))_{k=0}^{2^d-1}$ if and only if they
differ by coordinate reflections. Consequently, for $d\ge2$ there are
exactly $2^{(d-2)(d-3)/2}$ distinct profiles.
\end{introthm}

The proof uses finite Fourier analysis on $\Z/2^d\Z$. The Fourier
transform of the profile is the product $\prod_i(1+\chi_{w_i}+\chi_{-w_i})$
of characters, and the key point is that $1+t+t^{-1}$ does not vanish at
roots of unity of $2$-power order. For lattices tiling by discrete boxes
$\prod_i\{0,\ldots,n_i-1\}$ the factors become Dirichlet kernels, which
admit additional multiplicative symmetries. We determine them exactly.

\begin{introthm}[Theorem~\ref{thm:boxclass}]\label{thmC}
Let $d\ge2$, and let $n_1,\ldots,n_d\ge2$ have the same prime divisors.
The ordered cyclic profile determines every cocyclic upper triangular
lattice with diagonal $(n_1,\ldots,n_d)$ up to coordinate reflections if
and only if $n_1\in\{3,4,6\}$, or $n_1=2$ and either $d=2$ or $n_2=2$.
\end{introthm}

\subsection*{Organization}
Section~\ref{sec:partitions} introduces the intersection partitions and
their geometric interpretation. Section~\ref{sec:bases} treats triangular
bases, cocyclicity and the normalized congruence, and
Section~\ref{sec:reduction} gives the recursion in the dimension.
Sections~\ref{sec:signed} and~\ref{sec:binary} contain the connection with
Stern's sequence and the two rational trees, and
Section~\ref{sec:keplercubes} the geometric realization of the Kepler tree.
Sections~\ref{sec:profiles} and~\ref{sec:boxes} study cyclic profiles, for
cubes and for boxes. Section~\ref{sec:questions} collects open questions.
Appendix~\ref{app:groups} records the group-theoretic form of the
construction, and Appendix~\ref{app:infinite} a direct-limit version.

%======================================================================
\section{Intersection partitions of discrete cube tilings}
\label{sec:partitions}
%======================================================================
\subsection{Subcubes and a counting identity}
Write $\Delta_d=K_d-K_d=\{-1,0,1\}^d$. For $y\in\Delta_d$ set
\begin{equation}\label{eq:face}
 K_d(y)=K_d\cap(K_d-y)=\{u\in K_d:u+y\in K_d\}.
\end{equation}
Its $i$th factor is $\{1\}$ if $y_i=-1$, $\{0,1\}$ if $y_i=0$, and
$\{0\}$ if $y_i=1$. Thus $y\mapsto K_d(y)$ is a bijection between
$\Delta_d$ and the family $\mathscr K_d$ of all nonempty subcubes (faces)
of $K_d$; in particular $|\mathscr K_d|=3^d$. We shall often encode a face
by a word over $\{0,1,*\}$, where $*$ marks a free coordinate: $K_d(y)$
has the letter $1$, $*$ or $0$ in position $i$ according as $y_i=-1$, $0$
or $1$. For example, $K_3(1,0,-1)=\{0\}\times\{0,1\}\times\{1\}$ is the
edge $0{*}1$.

Let $T\subseteq\Z^d$ and suppose that $\{K_d+t:t\in T\}$ is a partition of
$\Z^d$. We refer to this as a \emph{discrete cube tiling} with translation
set $T$. For $x\in\Z^d$, define
\begin{align}
 P_x^T&=\{(K_d+x)\cap(K_d+t):t\in T\}\setminus\{\varnothing\},\label{eq:Px}\\
 \Pi_x^T&=P_x^T-x
   =\{K_d(x-t):t\in T,\ x-t\in\Delta_d\},\label{eq:Pix}\\
 \ell_T(x)&=|P_x^T|=|\Pi_x^T|
   =|T\cap(x+\Delta_d)|.\label{eq:counter}
\end{align}
Both $P_x^T$ and $\Pi_x^T$ are partitions, of $K_d+x$ and $K_d$,
respectively; translation of a family means translation of each member.
We call $\ell_T$ the \emph{partition counter} of the tiling.

\begin{theorem}\label{thm:sum}
For every $z\in\Z^d$, the families $\Pi_a^T$, $a\in K_d+z$, partition
$\mathscr K_d$. Consequently,
\begin{equation}\label{eq:sum}
 \sum_{a\in K_d+z}\ell_T(a)=3^d.
\end{equation}
\end{theorem}
\begin{proof}
Central symmetry gives $-K_d=K_d-\ones$, where $\ones=(1,\ldots,1)$.
Reflecting and translating the original tiling shows that
$\{K_d-t:t\in T\}$ also tiles $\Z^d$. For each $y\in\Delta_d$ there is
therefore a unique representation $y-z=b-t$ with $b\in K_d$ and $t\in T$.
Set $a=z+b$. Then $y=a-t$, so $K_d(y)\in\Pi_a^T$. Conversely, every
occurrence of $K_d(y)$ in some $\Pi_a^T$, $a\in K_d+z$, has this form, by
the injectivity of $y\mapsto K_d(y)$. This proves existence and uniqueness.
Counting the members of $\mathscr K_d$ gives \eqref{eq:sum}.
\end{proof}

\begin{corollary}\label{cor:mean}
Every translate $K_d+z$ contains a point $a$ with $\ell_T(a)\le(3/2)^d$
and a point $b$ with $\ell_T(b)\ge(3/2)^d$.
\end{corollary}

\subsection{Black cubes}\label{sec:black}
The partitions have a transparent geometric meaning. Place a closed unit
cube $B_t=t+[0,1]^d$ at every $t\in T$, and call these cubes
\emph{black}. The vertex set of $B_t$ is $K_d+t$, so every point of $\Z^d$
is a vertex of exactly one black cube. In particular, black cubes are
pairwise disjoint: a nonempty intersection of two lattice-aligned unit
cubes is a common face and contains a common vertex.

\begin{figure}[ht]
\centering
\begin{tikzpicture}[scale=0.82]
\draw[gray!45] (-1,-1) grid (6,5);
\foreach \x/\y in {0/0,2/0,4/0,-1/2,1/2,3/2,5/2,0/4,2/4,4/4}
 {\fill[black!78] (\x,\y) rectangle ++(1,1);
  \node[white,font=\scriptsize] at (\x+0.5,\y+0.5) {$1$};}
\foreach \y in {-1,1,3}
 {\foreach \x in {-1,...,5} {\node[font=\small] at (\x+0.5,\y+0.5) {$3$};}}
\foreach \x/\y in {-1/0,1/0,3/0,5/0,0/2,2/2,4/2,-1/4,1/4,3/4,5/4}
 {\node[font=\small] at (\x+0.5,\y+0.5) {$2$};}
\draw[red!75!black,very thick,dashed] (0,0) rectangle (2,2);
\draw[blue!75!black,line width=2.2pt] (3,3)--(4,3);
\fill[blue!75!black] (3,4) circle (3.2pt) (4,4) circle (3.2pt);
\draw[orange!90!black,line width=2.2pt] (3,0)--(3,1) (4,0)--(4,1);
\end{tikzpicture}
\caption{Black unit squares at the points of the lattice generated by
$(2,0)$ and $(1,2)$. Every integer point is a vertex of exactly one black
square. The number in the square $C_x$ is $\ell(x)$, the number of black
squares touching it. The white square $C_{(3,3)}$ is touched along its
lower edge and at its two upper corners (blue), giving the partition
$\{{*}0,01,11\}$ of its vertex set; $C_{(3,0)}$ is touched along its two
vertical edges (orange), giving $\{0{*},1{*}\}$. The dashed block of
the four squares $C_a$, $a\in K_2$ (its lower left corner is the origin),
illustrates Theorem~\ref{thm:sum}: $1+2+3+3=9$.}
\label{fig:black}
\end{figure}
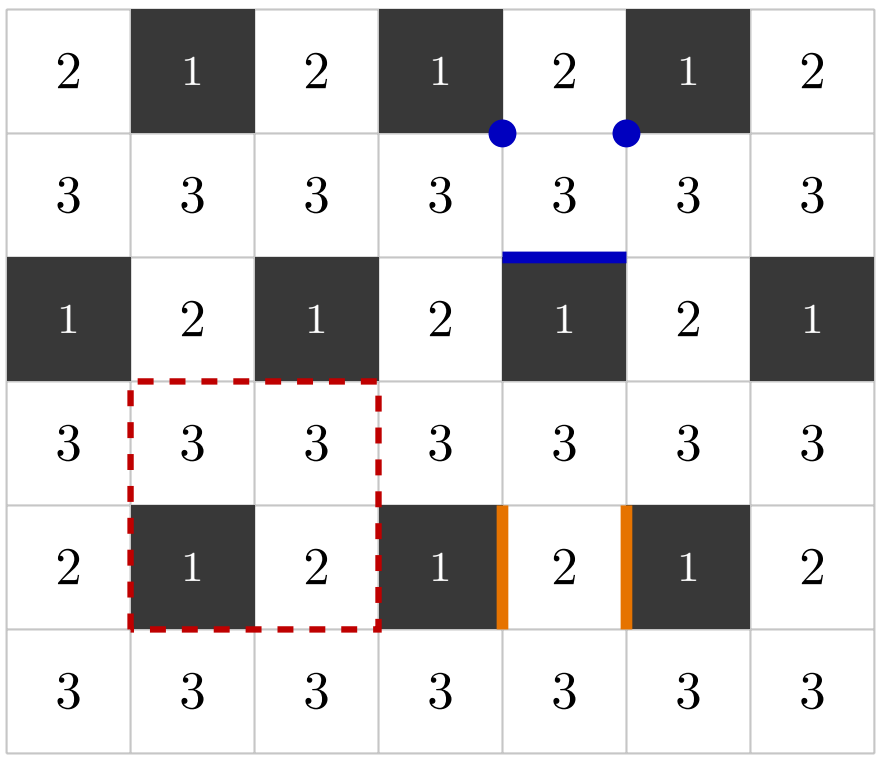

For $x\in\Z^d$ consider the unit cube $C_x=x+[0,1]^d$. Every nonempty
intersection $C_x\cap B_t$ is a face of $C_x$ with vertex set
$(K_d+x)\cap(K_d+t)$. Hence $P_x^T$ is the family of vertex sets of the
faces along which black cubes touch $C_x$, and
\[
 \ell_T(x)=\text{the number of black cubes touching }C_x.
\]
These faces partition the vertex set of $C_x$, though not necessarily the
whole cube $C_x$. If $x\in T$, then $C_x$ is black and $\ell_T(x)=1$; the
remaining unit cubes are \emph{white}. Figure~\ref{fig:black} shows a
planar example.

The same average appears in a containment count. For $x\in\Z^d$ let
$m_T(x)$ be the number of black cubes contained in $x+[0,3]^d$. A black
cube $B_t$ lies in $x+[0,3]^d$ if and only if $t-x\in\{0,1,2\}^d$, that
is, if and only if it touches $C_{x+\ones}$. Hence
\begin{equation}\label{eq:containment}
 m_T(x)=|T\cap(x+\{0,1,2\}^d)|=\ell_T(x+\ones),
 \qquad \sum_{a\in K_d+z}m_T(a)=3^d.
\end{equation}

\subsection{The lattice case}
Suppose now that $T=\Lambda$ is a full-rank sublattice of $\Z^d$. The
tiling condition is equivalent to $K_d$ being a complete system of
representatives for $\Z^d/\Lambda$. Denote the representative map by
$\rho_\Lambda:\Z^d\to K_d$, and transfer addition to $K_d$ by
$a\oplus b=\rho_\Lambda(a+b)$. We suppress the subscript $\Lambda$ when
there is no ambiguity. For $x\in\Z^d$ put
\begin{equation}\label{eq:displacements}
 L_x=\Delta_d\cap(x+\Lambda).
\end{equation}
Equations~\eqref{eq:Pix}--\eqref{eq:counter} become
\begin{equation}\label{eq:latticecounter}
 \Pi_x=\{K_d(y):y\in L_x\},\qquad
 \ell_\Lambda(x)=|L_x|.
\end{equation}

\begin{proposition}\label{prop:displacement}
For every $x\in\Z^d$,
\begin{equation}\label{eq:rho}
 L_x=\{\rho(x+u)-u:u\in K_d\}.
\end{equation}
On each part $K_d(y)\in\Pi_x$, the permutation $u\mapsto\rho(x+u)$ of
$K_d$ is the ordinary translation by $y$. Moreover,
\begin{equation}\label{eq:periodsym}
 \ell(x+\lambda)=\ell(x)\quad(\lambda\in\Lambda),
 \qquad \ell(-x)=\ell(x).
\end{equation}
\end{proposition}
\begin{proof}
Every displacement in \eqref{eq:rho} belongs to $\Delta_d\cap(x+\Lambda)$.
Conversely, write $y=v-u$ with $u,v\in K_d$ and $y\equiv x\pmod\Lambda$.
Then $v=\rho(x+u)$. If $u\in K_d(y)$, both $u+y$ and $\rho(x+u)$ belong to
$K_d$ and represent the same coset; hence they are equal. Periodicity
follows from \eqref{eq:displacements}, and negation is a bijection
$L_x\to L_{-x}$.
\end{proof}

\begin{proposition}\label{prop:product}
If $\Lambda=\Lambda_1\times\Lambda_2$ is a product of cube tiling
lattices, then $\ell_\Lambda(x_1,x_2)=\ell_{\Lambda_1}(x_1)\ell_{\Lambda_2}(x_2)$.
In particular,
$\ell_{2\Z^d}(x)=2^{|\{i:\,x_i\text{ is odd}\}|}$.
\end{proposition}
\begin{proof}
The displacement set \eqref{eq:displacements} factors as a Cartesian
product. In one dimension, the even coset contains only $0$ from
$\{-1,0,1\}$, while the odd coset contains both $-1$ and $1$.
\end{proof}

%======================================================================
\section{Triangular bases, cocyclicity and the normalized congruence}
\label{sec:bases}
%======================================================================
\subsection{Triangular bases}
The discrete tiling condition for a lattice $\Lambda$ is equivalent to
the half-open cubes $[0,2)^d+\lambda$, $\lambda\in\Lambda$, tiling
$\R^d$: each such cube is the disjoint union of the cubes
$[0,1)^d+u+\lambda$, $u\in K_d$. By the Minkowski--Haj\'os theorem
\cite{Hajos,SteinSzabo,Kolountzakis}, after a permutation of the coordinate
axes a cube tiling lattice has an upper triangular basis matrix with
diagonal entries $2$.

It costs nothing, and will be needed in Section~\ref{sec:boxes}, to work
with boxes from the start. Fix integers $n_1,\ldots,n_d\ge2$, and put
\[
 K(n)=\prod_{i=1}^d\{0,1,\ldots,n_i-1\},\qquad
 \Delta(n)=K(n)-K(n)=\prod_{i=1}^d\{-(n_i-1),\ldots,n_i-1\}.
\]
Thus $K_d=K(2,\ldots,2)$ and $\Delta_d=\Delta(2,\ldots,2)$. Matrices act
on column vectors, and the columns of a basis matrix $A$ generate the
lattice $A\Z^d$; we write $f_1,\ldots,f_d$ for the columns of $A$.

\begin{proposition}\label{prop:normalized}
Let $A=(a_{ij})$ be an upper triangular integer matrix with diagonal
$n_1,\ldots,n_d$.
\begin{enumerate}[(a)]
\item $K(n)$ is a complete system of representatives of $\Z^d/A\Z^d$.
In particular, $\{K(n)+\lambda:\lambda\in A\Z^d\}$ is a tiling of $\Z^d$,
$[\Z^d:A\Z^d]=n_1\cdots n_d$ and $A\Z^d\cap\Delta(n)=\{0\}$.
\item The lattice $A\Z^d$ has a unique upper triangular basis with the same
diagonal and $0\le a_{ij}<n_i$ for $i<j$.
\end{enumerate}
\end{proposition}
\begin{proof}
(a) Starting at $y=x\in\Z^d$, perform for $j=d,d-1,\ldots,1$ the update
\begin{equation}\label{eq:reduction}
 y\longleftarrow y-\left\lfloor\frac{y_j}{n_j}\right\rfloor f_j.
\end{equation}
The $j$th step makes $y_j\in\{0,\ldots,n_j-1\}$ and does not change the
coordinates $j+1,\ldots,d$. The result lies in $K(n)$ and is congruent to
$x$. If a nonzero vector $Az$ has last nonzero coefficient $z_j$, its
$j$th coordinate is $n_jz_j$, which lies outside
$\{-(n_j-1),\ldots,n_j-1\}$. Thus $A\Z^d\cap\Delta(n)=\{0\}$, proving the
uniqueness of representatives.

(b) For each $j$ and $i=j-1,j-2,\ldots,1$, subtract from $f_j$ a multiple
of $f_i$ so that its $i$th entry lies in $\{0,\ldots,n_i-1\}$; this does
not affect the rows below $i$. If two such matrices define the same
lattice, the difference of their $j$th columns lies in the lattice, has
zeros in the rows $i\ge j$ and entries in $\{-(n_i-1),\ldots,n_i-1\}$ in
the rows $i<j$. It lies in $\Delta(n)$, and hence vanishes by~(a).
\end{proof}

Let $\mathcal T_d$ be the set of upper triangular integer matrices with
diagonal entries $2$ and above-diagonal entries in $\{0,1\}$. A lattice
$A\Z^d$ with $A\in\mathcal T_d$ is called \emph{standard}. The coordinate
order is fixed in this terminology. By Proposition~\ref{prop:normalized},
$A\mapsto A\Z^d$ is injective on $\mathcal T_d$, so there are exactly
$2^{d(d-1)/2}$ standard cube tiling lattices.

\subsection{Cocyclic lattices and the normalized congruence}
A sublattice $\Lambda\subseteq\Z^d$ is \emph{cocyclic} if the quotient group
$\Z^d/\Lambda$ is cyclic. We say that $n_1,\ldots,n_d$ have \emph{common
prime support} if every prime dividing one of them divides all of them;
for instance, $n_i=2$ for all $i$. Throughout we write
\[
 N_i=n_1n_2\cdots n_i,\qquad N=N_d,\qquad M_i=N/N_i=n_{i+1}\cdots n_d,
\]
$\ZN=\Z/N\Z$, and $\nu_p(m)$ for the exponent of the prime $p$ in $m$.
The order of a residue $a\in\ZN$ is $\ord(a)=N/\gcd(a,N)$.

\begin{theorem}\label{thm:weights}
Let $A=(a_{ij})$ be an upper triangular integer matrix whose diagonal
entries $n_1,\ldots,n_d\ge2$ have common prime support, let
$\Lambda=A\Z^d$, and put $Q=\Z^d/\Lambda$. The following are equivalent:
\begin{enumerate}[(a)]
\item $Q$ is cyclic, that is, $\Lambda$ is cocyclic;
\item $e_d+\Lambda$ has order $N$ (equivalently, it generates $Q$);
\item $\gcd(a_{i,i+1},N)=1$ for $1\le i<d$.
\end{enumerate}
If they hold, there is a unique $w=(w_1,\ldots,w_d)\in\Z^d$ with
$0\le w_i<N$ and $w_d=1$ such that
\begin{equation}\label{eq:kernel}
 \Lambda=\Bigl\{x\in\Z^d:\ \sum_{i=1}^dw_ix_i\equiv0\pmod N\Bigr\}.
\end{equation}
Moreover, $\ord(w_i)=N_i$ in $\ZN$ for every $i$.
\end{theorem}

We call $w$ the \emph{normalized weights} of $\Lambda$.

\begin{proof}
By Proposition~\ref{prop:normalized}, $|Q|=N$.

(a)$\Rightarrow$(c). Suppose that a prime $p\mid N$ divides some
$a_{i,i+1}$. By the common prime support, $p$ divides every diagonal entry,
so $A$ is strictly upper triangular modulo $p$. All $(d-1)\times(d-1)$
minors of $A$ then vanish modulo $p$, except possibly the one obtained by
deleting the last row and the first column: every other minor retains the
zero first column or the zero last row. Modulo $p$, the remaining minor is
the determinant of an upper triangular matrix with diagonal
$a_{12},a_{23},\ldots,a_{d-1,d}$, so it also vanishes modulo $p$. Hence every
cofactor of $A$ is divisible by $p$. Since $\det A=N$,
\[
 \frac Np\,A^{-1}=\frac1p\operatorname{adj}A
\]
is an integer matrix. Therefore $(N/p)x\in\Lambda$ for every $x\in\Z^d$,
that is, every element of $Q$ has order dividing $N/p<N$. This contradicts
the cyclicity of the group $Q$ of order $N$. Hence (c) holds.

(c)$\Rightarrow$(b) and the weights. There is a unique $q\in\Q^d$ with
\begin{equation}\label{eq:qdef}
 q\cdot f_j=N\qquad(1\le j\le d),
\end{equation}
namely $q^{\mathsf T}=N(1,\ldots,1)A^{-1}$. Since $\det A=N$, we have
$NA^{-1}=\operatorname{adj}A$, so $q$ is integral. Geometrically, the
basis points $f_1,\ldots,f_d$ lie on the affine hyperplane $q\cdot x=N$.
We claim that for every prime $p\mid N$,
\begin{equation}\label{eq:valuation}
 \nu_p(q_i)=\nu_p(M_i)\qquad(1\le i\le d).
\end{equation}
The first equation of \eqref{eq:qdef} reads $n_1q_1=N$, so $q_1=M_1$.
Suppose that \eqref{eq:valuation} holds below $i\ge2$. The $i$th equation
reads
\begin{equation}\label{eq:qrec}
 n_iq_i=N-\sum_{j<i}a_{ji}q_j.
\end{equation}
The term $a_{i-1,i}q_{i-1}$ has $p$-adic valuation
$\nu_p(M_{i-1})=\nu_p(n_i)+\nu_p(M_i)$, because (c) makes $a_{i-1,i}$ prime
to $p$.
Every term with $j<i-1$ has valuation at least $\nu_p(M_j)$, and $N=M_0$.
For $0\le j<i-1$ we have $M_j=n_{j+1}\cdots n_{i-1}M_{i-1}$, a product with
at least one more factor divisible by $p$, so all these valuations exceed
$\nu_p(M_{i-1})$. Hence the right-hand
side of \eqref{eq:qrec} has valuation exactly $\nu_p(n_i)+\nu_p(M_i)$,
which proves \eqref{eq:valuation}. Consequently $\gcd(q_i,N)=M_i$; in
particular $q_d$ is prime to $N$.

Let $u\in\{0,\ldots,N-1\}$ satisfy $uq_d\equiv1\pmod N$, and let
$w_i\in\{0,\ldots,N-1\}$ be determined by $w_i\equiv uq_i\pmod N$. Then
$w_d=1$, and \eqref{eq:qdef} shows that every column of $A$ satisfies the
congruence in \eqref{eq:kernel}. The homomorphism
$x\mapsto\sum_iw_ix_i\bmod N$ from $\Z^d$ to $\ZN$ is onto, since
$w_d=1$, so its kernel has index $N$. This kernel contains $\Lambda$,
which also has index $N$; hence the two are equal. Thus the homomorphism
induces an isomorphism $Q\cong\ZN$ sending $e_d+\Lambda$ to $1$. This
proves (b). Multiplication by the unit $u$ does not change
$\gcd(q_i,N)=M_i$, so $\ord(w_i)=N/M_i=N_i$.

(b)$\Rightarrow$(a). This is immediate, since $|Q|=N$.

Uniqueness. If $v$ is another such tuple, then $e_i-w_ie_d\in\Lambda$ for
$i<d$, and the congruence defined by $v$ gives $v_i\equiv w_i\pmod N$.
Both lie in $\{0,\ldots,N-1\}$, so $v_i=w_i$.
\end{proof}

The converse construction shows that the normalized weights are a
complete set of coordinates for cocyclic triangular lattices.

\begin{proposition}\label{prop:converse}
Let $w\in\Z^d$ satisfy $w_d=1$ and $\ord(w_i)=N_i$ in $\ZN$ for all $i$.
Then $\Lambda_w=\{x:\sum_iw_ix_i\equiv0\pmod N\}$ is a cocyclic lattice
having an upper triangular basis with diagonal $n_1,\ldots,n_d$.
Consequently, if $n_1,\ldots,n_d$ have common prime support, the map
$\Lambda\mapsto w$ is a bijection between the cocyclic lattices with an
upper triangular basis with diagonal $n_1,\ldots,n_d$ and the tuples
$(w_1,\ldots,w_{d-1})$ with $\ord(w_i)=N_i$. Their number is
\begin{equation}\label{eq:countlattices}
 \prod_{i=1}^{d-1}\varphi(N_i),
\end{equation}
where $\varphi$ is Euler's function.
\end{proposition}
\begin{proof}
The lattice is cocyclic because $x\mapsto w\cdot x\bmod N$ is onto.
Let $Q_i\le\ZN$ be the image of $\Z^i\times\{0\}$, the subgroup generated
by $w_1,\ldots,w_i$. A cyclic group has a unique subgroup of each order,
and $\ord(w_j)$ divides $\ord(w_i)$ for $j\le i$; hence $Q_i=\langle w_i
\rangle$ and $|Q_i|=N_i$. Put $\Lambda_i=\Lambda_w\cap(\Z^i\times\{0\})$
and let $h_i\Z$ be the image of $\Lambda_i$ under the $i$th coordinate.
Then $\Lambda_{i-1}$ is the kernel of this coordinate map on $\Lambda_i$,
and therefore $[\Z^i:\Lambda_i]=h_i[\Z^{i-1}:\Lambda_{i-1}]$. Since
$\Z^i/\Lambda_i\cong Q_i$, we get $h_1\cdots h_i=N_i$, that is,
$h_i=n_i$. Choosing $\lambda_i\in\Lambda_i$ with $i$th coordinate $n_i$
gives an upper triangular basis $\lambda_1,\ldots,\lambda_d$ of
$\Lambda_w$ with diagonal $n_1,\ldots,n_d$.

Together with Theorem~\ref{thm:weights}, this gives the bijection. A
cyclic group of order $N$ contains exactly $\varphi(N_i)$ elements of
order $N_i$.
\end{proof}

\begin{remark}\label{rem:hnfcount}
Counting Hermite normal forms gives the same number. By
Theorem~\ref{thm:weights} and Proposition~\ref{prop:normalized}(b), the
superdiagonal entry $a_{i,i+1}$ ranges over the $\varphi(n_i)$ units
modulo $n_i$, and every other entry $a_{ij}$ over $n_i$ values, which
gives $\prod_{i<d}\varphi(n_i)n_i^{d-1-i}$. Under the common prime
support, $\varphi(m)/m$ is the same for all $m=n_i$ and $m=N_i$, and the
two counts agree.
\end{remark}

\begin{corollary}\label{cor:cyclicbinary}
A standard lattice $A\Z^d$, $A\in\mathcal T_d$, is cocyclic if and only if
every superdiagonal entry of $A$ equals $1$. There are $2^{(d-1)(d-2)/2}$
standard cocyclic cube tiling lattices in $\Z^d$, and their normalized
weights are exactly the tuples with $w_d=1$ and $\nu_2(w_i)=d-i$
for $i<d$.
\end{corollary}
\begin{proof}
Here $n_i=2$ and $N=2^d$. The only unit modulo $2$ is $1$. An element of
$\Z_{2^d}$ has order $2^i$ if and only if its $2$-adic valuation is
$d-i$, and there are $\varphi(2^i)=2^{i-1}$ of them; the product over
$i<d$ is $2^{(d-1)(d-2)/2}$.
\end{proof}

\begin{example}\label{ex:three}
For $d=3$ the two standard cocyclic lattices have bases
\[
 \begin{pmatrix}2&1&0\\0&2&1\\0&0&2\end{pmatrix},\qquad
 B_3=\begin{pmatrix}2&1&1\\0&2&1\\0&0&2\end{pmatrix}.
\]
For the first, the plane through the columns is $4x_1+2x_2+3x_3=8$, so
$q=(4,2,3)$. Since $3^{-1}\equiv3\pmod 8$, we get $w=(4,6,1)$. For $B_3$
we get $q=(4,2,1)$ and $w=(4,2,1)$. The two weight vectors differ by the
sign of $w_2$; see Section~\ref{sec:profiles}.
\end{example}

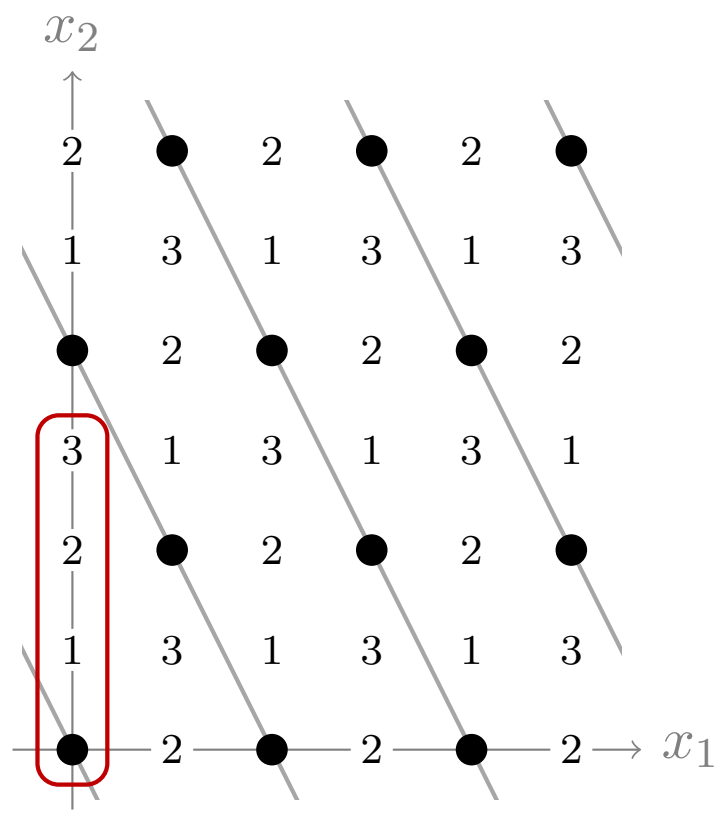
\begin{figure}[ht]
\centering
\begin{tikzpicture}[scale=0.66]
\draw[->,gray] (-0.6,0)--(5.7,0) node[right] {$x_1$};
\draw[->,gray] (0,-0.6)--(0,6.8) node[above] {$x_2$};
\begin{scope}
\clip (-0.5,-0.5) rectangle (5.5,6.5);
\foreach \k in {0,...,4} {\draw[thick,gray!70] (-1,4*\k+2)--(6,4*\k-12);}
\end{scope}
\foreach \x in {0,...,5} {\foreach \y in {0,...,6} {
 \pgfmathtruncatemacro{\r}{mod(2*\x+\y,4)}
 \ifnum\r=0 \fill (\x,\y) circle (4.5pt);
 \else \fill[white] (\x,\y) circle (6pt);
   \node[font=\scriptsize] at (\x,\y) {$\r$};
 \fi}}
\draw[red!75!black,thick,rounded corners] (-0.35,-0.35) rectangle (0.35,3.35);
\end{tikzpicture}
\caption{The lattice $\Lambda_2^b=\{x\in\Z^2:2x_1+x_2\equiv0\pmod 4\}$
(black dots) lies on the lines $2x_1+x_2=4k$. Every other integer point is
labeled by the residue of $2x_1+x_2$ modulo $4$, which identifies its
coset. The framed points $0,e_2,2e_2,3e_2$ represent the four cosets.}
\label{fig:hyperplanes}
\end{figure}

\begin{remark}[Parallel hyperplanes]\label{rem:hyperplanes}
Let $w$ be the normalized weights of a cocyclic lattice $\Lambda$, and let
$\mathcal H_m=\{x\in\R^d:w\cdot x=m\}$ for $m\in\Z$. Since $w_d=1$, the
integer points of $\mathcal H_m$ are
\[
 \mathcal H_m\cap\Z^d=\Bigl\{me_d+\sum_{i=1}^{d-1}t_i(e_i-w_ie_d):
 t_i\in\Z\Bigr\}.
\]
The coset of an integer point is determined by the index of the
hyperplane containing it: $x\in\mathcal H_m$ lies in $ke_d+\Lambda$ if and
only if $m\equiv k\pmod N$. In particular,
\[
 \Lambda=\bigcup_{k\in\Z}\bigl(\mathcal H_{kN}\cap\Z^d\bigr)
 =\langle e_1-w_1e_d,\ldots,e_{d-1}-w_{d-1}e_d,\ Ne_d\rangle_\Z,
\]
and the points $0,e_d,\ldots,(N-1)e_d$ represent all cosets. The
hyperplane $q\cdot x=N$ of the proof of Theorem~\ref{thm:weights} serves
only to construct $w$ and should not be confused with this family. See
Figure~\ref{fig:hyperplanes}.
\end{remark}

%======================================================================
\section{Dimension reduction}\label{sec:reduction}
%======================================================================
Let $A_d$ be any upper triangular integer matrix with diagonal entries
$2$, not necessarily normalized, and extend it by
\[
 A_{d+1}=\begin{pmatrix}A_d&g\\0&2\end{pmatrix},\qquad g\in\Z^d.
\]
Write $\Lambda_j=A_j\Z^j$ and $\ell_j=\ell_{\Lambda_j}$. Both lattices
tile by Proposition~\ref{prop:normalized}.

\begin{theorem}\label{thm:recursion}
For all $x\in\Z^d$ and $r\in\Z$,
\begin{align}
 \ell_{d+1}(x,2r)&=\ell_d(x-rg),\label{eq:even}\\
 \ell_{d+1}(x,2r+1)&=\ell_d(x-rg)+\ell_d(x-(r+1)g).\label{eq:odd}
\end{align}
In particular, $\ell_{d+1}(x,0)=\ell_d(x)$.
\end{theorem}
\begin{proof}
We count the vectors $(y,t)\in\Delta_d\times\{-1,0,1\}$ congruent to
$(x,s)$ modulo $\Lambda_{d+1}$. The vectors of $\Lambda_{d+1}$ have the
form $(A_dz+qg,2q)$. Hence $s-t=2q$, and $y\equiv x-qg\pmod{\Lambda_d}$.
If $s=2r$, then $t=0$ and $q=r$. If $s=2r+1$, there are exactly two
possibilities, $(t,q)=(1,r)$ and $(-1,r+1)$, and their counts add.
\end{proof}

Repeated application shows that appending zero coordinates does not
change the counter, provided that the smaller basis is the corresponding
leading principal submatrix of the larger one.
Equations~\eqref{eq:even}--\eqref{eq:odd} also give a recursive algorithm
for evaluating the counter without enumerating $\Delta_d$.

%======================================================================
\section{Signed bidiagonal lattices and Stern's sequence}\label{sec:signed}
%======================================================================
\subsection{Stern rows as cyclic profiles}
Fix signs $\varepsilon_i\in\{-1,1\}$, and let
\begin{equation}\label{eq:signed}
 \Lambda_d^{\varepsilon}
 =\langle 2e_1,\ 2e_2+\varepsilon_1e_1,\ldots,
                 2e_d+\varepsilon_{d-1}e_{d-1}\rangle_{\Z}.
\end{equation}
By Theorem~\ref{thm:weights}, these are cocyclic cube tiling lattices with
generator $e_d$. Set $c_d(k)=\ell_{\Lambda_d^{\varepsilon}}(ke_d)$, and
let $c_0\equiv1$.

\begin{theorem}\label{thm:stern}
The function $c_d$ is independent of the signs and has period $2^d$. It is
determined by $c_1(0)=1$, $c_1(1)=2$, and
\begin{equation}\label{eq:cyclicstern}
 c_{d+1}(2r)=c_d(r),\qquad
 c_{d+1}(2r+1)=c_d(r)+c_d(r+1)\qquad(r\in\Z).
\end{equation}
If $s$ is Stern's diatomic sequence (A002487 in \cite{OEIS}), $s(0)=0$,
$s(1)=1$ and
\begin{equation}\label{eq:stern}
 s(2r)=s(r),\qquad s(2r+1)=s(r)+s(r+1)\qquad(r\ge1),
\end{equation}
then
\begin{equation}\label{eq:sternrow}
 c_d(k)=s(2^d+k)\qquad(0\le k<2^d).
\end{equation}
\end{theorem}
\begin{proof}
Apply Theorem~\ref{thm:recursion} with $x=0$ and $g=\varepsilon_de_d$.
The evenness \eqref{eq:periodsym} of $\ell_d$ eliminates both the minus
sign and $\varepsilon_d$, giving \eqref{eq:cyclicstern}. The initial
values follow by direct counting in $\{-1,0,1\}$. This proves the
independence of the signs, and the order $2^d$ of $e_d$ gives the period.
For \eqref{eq:sternrow}, induction uses \eqref{eq:stern} at each even or
odd index. At the right-hand endpoint one uses $s(2^{d+1})=1=c_d(0)$, so
the periodic wrap agrees with the next endpoint of Stern's row.
\end{proof}

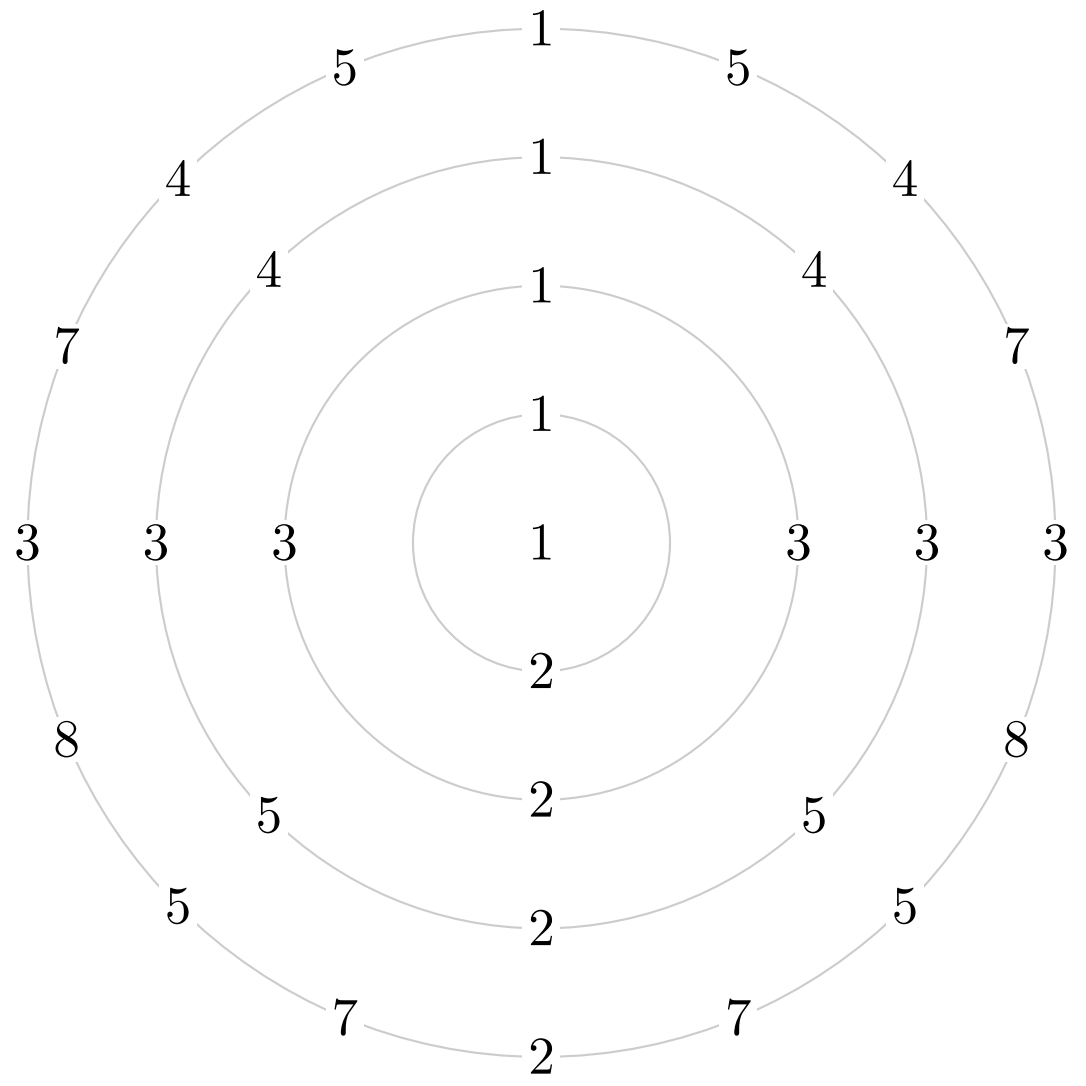
\begin{figure}[ht]
\centering
\begin{tikzpicture}
\foreach \r in {1,2,3,4} {\draw[gray!40] (0,0) circle ({0.85*\r});}
\node[font=\small] at (0,0) {$1$};
\foreach \v [count=\k from 0] in {1,2}
 {\node[fill=white,inner sep=1pt,font=\small] at ({90-180*\k}:0.85) {$\v$};}
\foreach \v [count=\k from 0] in {1,3,2,3}
 {\node[fill=white,inner sep=1pt,font=\small] at ({90-90*\k}:1.7) {$\v$};}
\foreach \v [count=\k from 0] in {1,4,3,5,2,5,3,4}
 {\node[fill=white,inner sep=1pt,font=\small] at ({90-45*\k}:2.55) {$\v$};}
\foreach \v [count=\k from 0] in {1,5,4,7,3,8,5,7,2,7,5,8,3,7,4,5}
 {\node[fill=white,inner sep=1pt,font=\small] at ({90-22.5*\k}:3.4) {$\v$};}
\end{tikzpicture}
\caption{The cyclic rows $c_0,c_1,\ldots,c_4$ on concentric circles, read
clockwise from the top. Each entry of $c_d$ reappears in $c_{d+1}$ at the
same angle, and the sum of two cyclically adjacent entries is inserted
between them. Rotating a ring corresponds to translating the labeling;
cf.\ Section~\ref{sec:keplercubes}.}
\label{fig:rings}
\end{figure}

The first rows are
\[
\begin{array}{c|l}
 d&c_d(0),c_d(1),\ldots,c_d(2^d-1)\\\hline
 1&1,2\\
 2&1,3,2,3\\
 3&1,4,3,5,2,5,3,4\\
 4&1,5,4,7,3,8,5,7,2,7,5,8,3,7,4,5.
\end{array}
\]
They have a cyclic insertion description: each entry is retained, and the
sum of two cyclically adjacent entries is inserted between them; see
Figure~\ref{fig:rings}.

\begin{remark}\label{rem:signedweights}
The normalized weights of $\Lambda_d^\varepsilon$ are easy to compute. The
generator $2e_j+\varepsilon_{j-1}e_{j-1}$ gives
$2w_j+\varepsilon_{j-1}w_{j-1}\equiv0$, so
$w_{j-1}\equiv-2\varepsilon_{j-1}w_j$ and
\[
 w_i\equiv(-2)^{d-i}\varepsilon_i\varepsilon_{i+1}\cdots\varepsilon_{d-1}
 \equiv\pm2^{d-i}\pmod{2^d}.
\]
As in Section~\ref{sec:profiles} below, $c_d(k)$ is the number of
$y\in\Delta_d$ with $\sum_iw_iy_i\equiv k$. Replacing $y_i$ by $-y_i$ where
necessary, $c_d(k)$ is the number of representations
$k\equiv\sum_{i=1}^dy_i2^{d-i}\pmod{2^d}$ with digits $y_i\in\{-1,0,1\}$.
This gives a second proof of the independence of the signs.
\end{remark}

\subsection{Calkin--Wilf labels on the representative cubes}
Let $\rho_d$ be the representative map of $\Lambda_d^\varepsilon$, for a
fixed infinite sequence of signs. Consider the rooted binary tree whose
level $d$ is $K_d$, with an extra empty root at level zero. The empty root
has children $0$ and $1$, and the children of $x\in K_d$, $d\ge1$, are
\begin{equation}\label{eq:signedchildren}
 (\rho_d(-\varepsilon_dx),0),\qquad (\rho_d(-\varepsilon_dx),1).
\end{equation}
Since multiplication by $-\varepsilon_d$ permutes the quotient, every
vertex of level $d+1$ has exactly one parent.

\begin{theorem}\label{thm:cw}
Label the empty root by $(1,1)$ and each $x\in K_d$ by
\begin{equation}\label{eq:pair}
 (\ell_d(x),\ell_d(x+e_d)),\qquad \ell_d=\ell_{\Lambda_d^\varepsilon}.
\end{equation}
If a vertex has label $(p,q)$, its children in \eqref{eq:signedchildren}
have labels $(p,p+q)$ and $(p+q,q)$. Replacing each label by the ratio of
its coordinates gives the Calkin--Wilf tree. In particular, every positive
rational number occurs exactly once.
\end{theorem}
\begin{proof}
The first level has labels $(1,2)$ and $(2,1)$. For $d\ge1$, put
$z=\rho_d(-\varepsilon_dx)$. By periodicity and evenness,
$\ell_d(z)=\ell_d(x)=p$ and $\ell_d(z-\varepsilon_de_d)=\ell_d(x+e_d)=q$.
Theorem~\ref{thm:recursion}, with $g=\varepsilon_de_d$, gives
\[
 \ell_{d+1}(z,0)=p,\qquad \ell_{d+1}(z,1)=p+q,\qquad \ell_{d+1}(z,2)=q,
\]
and these are the required labels. The rules preserve coprimality. A
coprime pair $(p,q)\ne(1,1)$ has a unique possible parent, $(p,q-p)$ if
$p<q$ and $(p-q,q)$ if $p>q$. Repeated passage to the parent decreases
$p+q$ and ends at $(1,1)$, which proves existence and uniqueness.
\end{proof}

\subsection{Normalized bases}
The first sign in \eqref{eq:signed} does not change the lattice, since
$2e_2-e_1$ and $2e_2+e_1$ differ by $2e_1$. The remaining signs have a
simple interpretation in the normalized basis. Write its $j$th column as
$2e_j+h_{j-1}$, with $h_0=0$. For $2\le j\le d$, reduction of the signed
column gives
\begin{equation}\label{eq:gradual}
 h_{j-1}=\begin{cases}
 e_{j-1},&\varepsilon_{j-1}=1,\\
 h_{j-2}+e_{j-1},&\varepsilon_{j-1}=-1;
 \end{cases}
\end{equation}
in the second case one adds the preceding normalized column. Thus the
block of consecutive ones in a normalized column either starts anew at
coordinate $j-1$ or extends the block of the preceding column. All signs
equal to $-1$ give the matrix $B_d$ of Section~\ref{sec:binary}, with all
above-diagonal entries equal to $1$.

\begin{corollary}\label{cor:mixedlattices}
For $d\ge2$, the construction \eqref{eq:gradual} gives $2^{d-2}$ distinct
standard lattices. All of them have the cyclic profile
\eqref{eq:sternrow} along $e_d$, and hence the same multiset of counter
values on $K_d$.
\end{corollary}
\begin{proof}
For $j\ge3$ the two choices in \eqref{eq:gradual} are distinct, since
$h_{j-2}$ contains $e_{j-2}$. Normalized bases determine distinct lattices
by Proposition~\ref{prop:normalized}. Theorem~\ref{thm:stern} gives the
common sequence, and $e_d$ generates the quotient.
\end{proof}

These $2^{d-2}$ lattices form a single orbit under coordinate reflections
(Remark~\ref{rem:sternorbit}), and by Theorem~\ref{thm:main} no other
standard lattice has the Stern profile.

%======================================================================
\section{Binary words, the Calkin--Wilf tree and the Kepler tree}
\label{sec:binary}
%======================================================================
Let $B_d$ be the upper triangular $d\times d$ matrix with diagonal entries
$2$ and all entries above the diagonal equal to $1$, put
$\Lambda_d^b=B_d\Z^d$, and
\[
 F_d(x)=\sum_{i=1}^d 2^{d-i}x_i.
\]
Each column of $B_d$ is mapped by $F_d$ to $2^d$, so
$\Lambda_d^b\subseteq\ker(F_d\bmod 2^d)$. Both lattices have index $2^d$,
so they are equal; that is, $F_d$ gives the normalized weights
$(2^{d-1},\ldots,2,1)$ of $\Lambda_d^b$. The induced group operation on
the binary words in $K_d$ is ordinary binary addition modulo $2^d$, with
the leftmost bit most significant. For a nonzero word whose rightmost $1$
is in position $j$, its inverse is
\[
 (a_1,\ldots,a_{j-1},1,0,\ldots,0)'=(1-a_1,\ldots,1-a_{j-1},1,0,\ldots,0),
\]
which must not be confused with the full bitwise complement. Since
$\Lambda_d^b$ is one of the signed lattices of Section~\ref{sec:signed},
\begin{equation}\label{eq:binaryprofile}
 \ell_{\Lambda_d^b}(x)=c_d(F_d(x))\qquad(x\in\Z^d).
\end{equation}

For a binary word $v$ write $\bar v$ for its bitwise complement and
$\pi(v)=\ell_{\Lambda_{|v|}^b}(v)$; set $\pi(\varnothing)=1$.
Juxtaposition denotes concatenation. The last column of $B_{d+1}$ has
$g=\ones$, and $v-\ones=-\bar v$, so Theorem~\ref{thm:recursion} and
evenness yield
\begin{equation}\label{eq:binaryrec}
 \pi(v0)=\pi(v),\qquad \pi(v1)=\pi(v)+\pi(\bar v).
\end{equation}
In particular $\pi(0)=1$, $\pi(1)=2$, and appending zeros leaves $\pi$
unchanged. For example,
$\pi(10101)=\pi(1010)+\pi(0101)=5+8=13$.

\begin{proposition}\label{prop:binarycw}
On the complete binary tree of words, with children $v0$ and $v1$, the
labels $(\pi(v),\pi(\bar v))$ satisfy the Calkin--Wilf rules. Hence
$v\mapsto\pi(v)/\pi(\bar v)$ is the Calkin--Wilf tree, the word $v$
describing the path from the root ($0$ = left, $1$ = right).
\end{proposition}
\begin{proof}
The root has label $(1,1)$. If $v$ has label $(p,q)$, then
$\overline{v0}=\bar v1$ and $\overline{v1}=\bar v0$, and
\eqref{eq:binaryrec} gives the labels $(p,p+q)$ and $(p+q,q)$.
\end{proof}

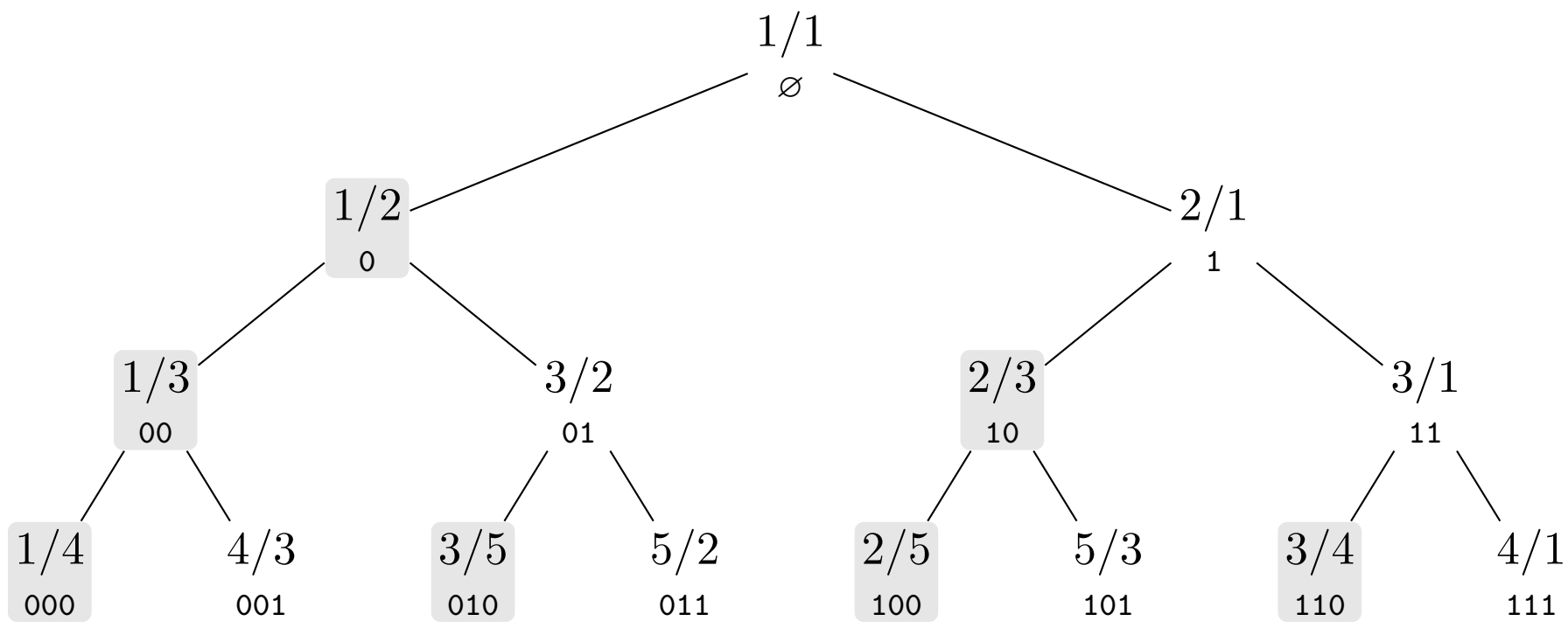
\begin{figure}[ht]
\centering
\begin{tikzpicture}[level distance=13mm,
 level 1/.style={sibling distance=64mm},
 level 2/.style={sibling distance=32mm},
 level 3/.style={sibling distance=16mm},
 every node/.style={inner sep=1.5pt,align=center,font=\small},
 kep/.style={fill=gray!20,rounded corners=2pt}]
\node {$1/1$\\[-1pt]\tiny$\varnothing$}
 child {node[kep] {$1/2$\\[-1pt]\tiny\texttt{0}}
  child {node[kep] {$1/3$\\[-1pt]\tiny\texttt{00}}
   child {node[kep] {$1/4$\\[-1pt]\tiny\texttt{000}}}
   child {node {$4/3$\\[-1pt]\tiny\texttt{001}}}}
  child {node {$3/2$\\[-1pt]\tiny\texttt{01}}
   child {node[kep] {$3/5$\\[-1pt]\tiny\texttt{010}}}
   child {node {$5/2$\\[-1pt]\tiny\texttt{011}}}}}
 child {node {$2/1$\\[-1pt]\tiny\texttt{1}}
  child {node[kep] {$2/3$\\[-1pt]\tiny\texttt{10}}
   child {node[kep] {$2/5$\\[-1pt]\tiny\texttt{100}}}
   child {node {$5/3$\\[-1pt]\tiny\texttt{101}}}}
  child {node {$3/1$\\[-1pt]\tiny\texttt{11}}
   child {node[kep] {$3/4$\\[-1pt]\tiny\texttt{110}}}
   child {node {$4/1$\\[-1pt]\tiny\texttt{111}}}}};
\end{tikzpicture}
\caption{The Calkin--Wilf tree as $v\mapsto\pi(v)/\pi(\bar v)$
(Proposition~\ref{prop:binarycw}); the word $v$ is shown below each
fraction. It is also the tree of ratios of the labels in
Theorem~\ref{thm:cw}. The shaded vertices, the words ending in $0$, carry
the fractions of the Kepler tree (Remark~\ref{rem:keplerincw}).}
\label{fig:cw}
\end{figure}

\begin{theorem}\label{thm:kepler}
Take the nonempty binary words ending in $0$, root the tree at $0$, and
give $v$ the children $v0$ and $\bar v0$. The labels
\begin{equation}\label{eq:keplerlabel}
 \frac{\pi(v)}{\pi(\bar v)}
\end{equation}
form the Kepler tree rooted at $1/2$. They enumerate $\Q\cap(0,1)$ without
repetition.
\end{theorem}
\begin{proof}
Every word of length at least two ending in $0$ has a unique parent:
remove its final zero, and complement the remaining word if that word ends
in $1$. If the pair at $v$ is $(p,q)$, then \eqref{eq:binaryrec} gives the
pairs $(p,p+q)$ at $v0$ and $(q,p+q)$ at $\bar v0$. The root has pair
$(1,2)$. These rules preserve $0<p<q$ and coprimality. Conversely, a
reduced fraction $p/q\in(0,1)$ other than $1/2$ has the unique parent
$p/(q-p)$ if $2p<q$, and $(q-p)/p$ if $2p>q$. The denominator strictly
decreases, so the descent terminates at $1/2$.
\end{proof}

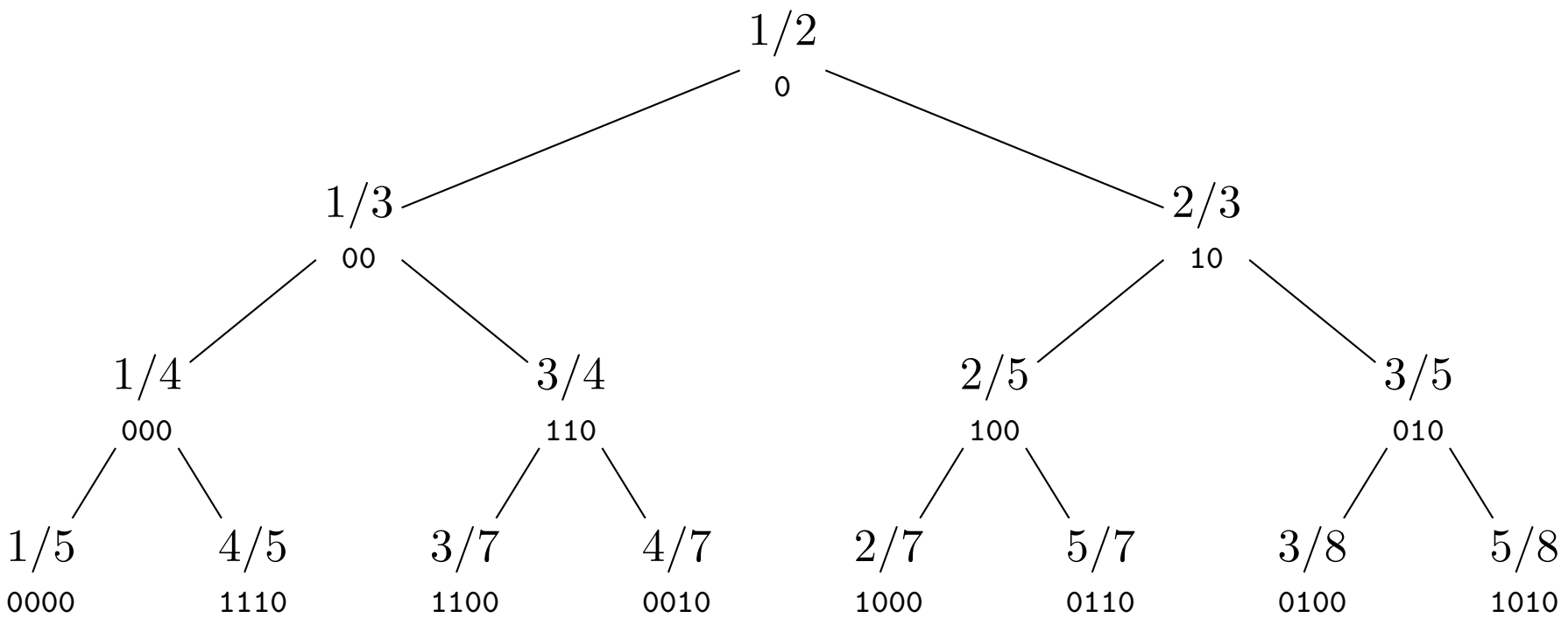
\begin{figure}[ht]
\centering
\begin{tikzpicture}[level distance=13mm,
 level 1/.style={sibling distance=64mm},
 level 2/.style={sibling distance=32mm},
 level 3/.style={sibling distance=16mm},
 every node/.style={inner sep=1.5pt,align=center,font=\small}]
\node {$1/2$\\[-1pt]\tiny\texttt{0}}
 child {node {$1/3$\\[-1pt]\tiny\texttt{00}}
  child {node {$1/4$\\[-1pt]\tiny\texttt{000}}
   child {node {$1/5$\\[-1pt]\tiny\texttt{0000}}}
   child {node {$4/5$\\[-1pt]\tiny\texttt{1110}}}}
  child {node {$3/4$\\[-1pt]\tiny\texttt{110}}
   child {node {$3/7$\\[-1pt]\tiny\texttt{1100}}}
   child {node {$4/7$\\[-1pt]\tiny\texttt{0010}}}}}
 child {node {$2/3$\\[-1pt]\tiny\texttt{10}}
  child {node {$2/5$\\[-1pt]\tiny\texttt{100}}
   child {node {$2/7$\\[-1pt]\tiny\texttt{1000}}}
   child {node {$5/7$\\[-1pt]\tiny\texttt{0110}}}}
  child {node {$3/5$\\[-1pt]\tiny\texttt{010}}
   child {node {$3/8$\\[-1pt]\tiny\texttt{0100}}}
   child {node {$5/8$\\[-1pt]\tiny\texttt{1010}}}}};
\end{tikzpicture}
\caption{The Kepler tree realized by the words ending in $0$
(Theorem~\ref{thm:kepler}); the children of $v$ are $v0$ and $\bar v0$.}
\label{fig:kepler}
\end{figure}

\begin{remark}\label{rem:keplerincw}
The Kepler label of a word $v$ ending in $0$ coincides with its
Calkin--Wilf label from Proposition~\ref{prop:binarycw}. Thus the Kepler
fractions are the left children in the Calkin--Wilf tree, and the Kepler
children of $v$ are the left children of $v$ and of its mirror image
$\bar v$, whose Calkin--Wilf label is the reciprocal $\pi(\bar v)/\pi(v)$.
Compare Figures~\ref{fig:cw} and~\ref{fig:kepler}.
\end{remark}

\subsection{Mixing two rational-tree rules}
The two branching rules can also be mixed level by level. Start at
$(1,1)$, and choose independently at each level one of the rules
\[
 \beta_1(p,q)=\bigl((p,p+q),(p+q,q)\bigr),\qquad
 \beta_2(p,q)=\bigl((q,p+q),(p+q,p)\bigr),
\]
which is then applied to every vertex of that level.

\begin{proposition}\label{prop:mixed}
Every such choice produces, at each level, exactly the same multiset of
ordered pairs as the Calkin--Wilf tree.
\end{proposition}
\begin{proof}
Let $\sigma(p,q)=(q,p)$. The multiset at a Calkin--Wilf level is
$\sigma$-invariant, by induction: swapping the parent coordinates swaps
the two child coordinates and their order. Moreover,
$\beta_2=\beta_1\circ\sigma$. Therefore $\beta_1$ and $\beta_2$ have
identical multiset images on every $\sigma$-invariant level, and induction
proves the assertion.
\end{proof}

The proof requires the same rule at all vertices of a level; it does not
cover an independent choice of rule at each vertex.

%======================================================================
\section{Kepler cubes and generating face partitions}\label{sec:keplercubes}
%======================================================================
In Theorem~\ref{thm:kepler} the Kepler tree came from the arithmetic of
the counter. In this section we show that it is also encoded in the
geometry of the labeled cube. We formulate local conditions on a labeling
of the vertices of $K_d$ (Definition~\ref{def:kepler}), a sum rule along
the diagonals of a partition of a face into subfaces and a matching between
these diagonals and edges, which make no reference to rational trees, and
show that they force the levels of the Kepler tree
(Proposition~\ref{prop:treeconsequence}). For the binary lattice, all
translated counter labelings satisfy these conditions, and the required
partitions are supplied by the tiling itself: they are partitions
$\Pi_x$ of Section~\ref{sec:partitions}
(Theorem~\ref{thm:allcubes}). Finally, every such partition of a white cube
is laminated in a unique direction (Theorem~\ref{thm:lamination}). The
product lattice $2\Z^2$ shows that the Kepler property is not automatic
for cube tilings (Example~\ref{ex:product}).

Throughout this section, $\Lambda_d^b$ and $F_d$ are as in
Section~\ref{sec:binary}, $\ell_d=\ell_{\Lambda_d^b}$, and $c_d$ is the
cyclic counter of Theorem~\ref{thm:stern}.

\subsection{Faces, diagonals and edge spectra}
We encode faces of $K_d$ by words in $\{0,1,*\}^d$ as in
Section~\ref{sec:partitions}. The \emph{antipodal map} $\tau_S$ of a face
$S$ complements all free coordinates and fixes the others. A
\emph{diagonal} of $S$ is an unordered pair $\{u,\tau_Su\}$ of distinct
vertices of $S$; a face of dimension $k\ge1$ has $2^{k-1}$ diagonals, and
a vertex has none. For a partition $\mathcal G$ of a face into
positive-dimensional faces, the diagonals of all parts partition the
vertex set of that face. Let $\sigma_k$ flip coordinate $k$.

For a labeling $a:K_d\to\Z_{>0}$ and a face $H$ whose first $j$
coordinates are free, let $\Ed_j(a;H)$ be the multiset of unordered
endpoint-label pairs $\{a(u),a(\sigma_ju)\}$ of the edges of $H$ parallel
to $e_j$, and write $\Ed_d(a)=\Ed_d(a;K_d)$. For a partition
$\mathcal G$ as above, let $\Dg(a;\mathcal G)$ be the multiset of the pairs
$\{a(u),a(v)\}$ over all diagonals $\{u,v\}$ of all parts. Each edge or
diagonal is counted once, and its label pair retains both entries, also
when they coincide.

\subsection{Kepler cubes}
\begin{definition}\label{def:kepler}
Let $d\ge1$. A \emph{Kepler cube} is a labeling $a:K_d\to\Z_{>0}$ admitting
a nested flag of faces
\[
 H_1\subset\cdots\subset H_d=K_d,\qquad
 H_k=\{0,1\}^k\times\{b_{k+1}\}\times\cdots\times\{b_d\},
\]
with fixed bits $b_2,\ldots,b_d\in\{0,1\}$, such that the endpoint labels
of the edge $H_1$ are $1$ and $2$, and such that for each $2\le k\le d$
there is a partition $\mathcal G_k$ of $H_{k-1}$ into positive-dimensional
faces with the following two properties.
\begin{enumerate}[(i)]
\item \emph{Sum rule.} For every diagonal $\{u,v\}$ of every part,
\begin{equation}\label{eq:generator}
 a(\sigma_ku)=a(\sigma_kv)=a(u)+a(v).
\end{equation}
\item \emph{Compatibility.} The diagonal labels match the edges of the
preceding direction:
\begin{equation}\label{eq:compatible}
 \Dg(a;\mathcal G_k)=\Ed_{k-1}(a;H_{k-1}).
\end{equation}
\end{enumerate}
We call $\mathcal G_k$ a \emph{generating partition}, and for $d\ge2$ we
call $\mathcal G_d$ a \emph{top generating partition}.
\end{definition}

Condition \eqref{eq:compatible} says that the diagonals of $\mathcal G_k$
can be matched bijectively with the $(k-1)$-edges of $H_{k-1}$, preserving
endpoint labels. The vertex $\sigma_ku$ lies in the facet of $H_k$
opposite to $H_{k-1}$, so the sum rule determines the labels of that facet.
No rational-tree condition is included in the definition, and the
normalization $\{1,2\}$ fixes the root $1/2$.

\begin{proposition}\label{prop:treeconsequence}
In a Kepler cube, the unordered endpoint-label pairs on the $k$-edges of
$H_k$ are exactly the numerator--denominator pairs at level $k-1$ of the
Kepler tree, each once; the root $1/2$ has level zero. In particular,
dividing the smaller endpoint label by the larger one gives that level.
\end{proposition}
\begin{proof}
The case $k=1$ is the normalization. Suppose the assertion holds for $k-1$.
By \eqref{eq:compatible}, each pair $\{p,q\}$ occurring on the
$(k-1)$-edges of $H_{k-1}$ occurs with the same multiplicity on the
diagonals of $\mathcal G_k$. By \eqref{eq:generator}, a diagonal
$\{u,v\}$ with labels $p,q$ produces the two $k$-edges
$\{u,\sigma_ku\}$ and $\{v,\sigma_kv\}$, with label pairs
\begin{equation}\label{eq:children}
 \{p,p+q\},\qquad\{q,p+q\}.
\end{equation}
For $p<q$ these are precisely the children of $p/q$ in
Theorem~\ref{thm:kepler}. The diagonals partition the vertices of
$H_{k-1}$, so every $k$-edge of $H_k$ is obtained exactly once.
\end{proof}

\subsection{The cyclic labelings and black partitions}
Define the cyclic family of labelings
\begin{equation}\label{eq:C-labels}
 a_d^t(u)=c_d(t+F_d(u)),\qquad u\in K_d,\quad t\in\Z/2^d\Z.
\end{equation}
By \eqref{eq:binaryprofile} and $F_d(\ones)=2^d-1\equiv-1$, the
translated counter labelings and the containment labelings
\eqref{eq:containment} belong to this family:
\begin{equation}\label{eq:translations}
 \ell_d(x+u)=a_d^{F_d(x)}(u),\qquad
 m_d(x+u)=\ell_d(x+u+\ones)=a_d^{F_d(x)-1}(u),
\end{equation}
where $m_d=m_{\Lambda_d^b}$. For a labeling of $K_{d+1}$ write $(A,B)$ for
its restrictions to the facets with last coordinate $0$ and $1$; addition
of labelings is pointwise. Since $F_{d+1}(u,b)=2F_d(u)+b$,
\eqref{eq:cyclicstern} gives
\begin{equation}\label{eq:cubrec}
 a_{d+1}^{2n}=\bigl(a_d^n,\ a_d^n+a_d^{n+1}\bigr),\qquad
 a_{d+1}^{2n+1}=\bigl(a_d^n+a_d^{n+1},\ a_d^{n+1}\bigr),
\end{equation}
with $a_1^0=(1,2)$ and $a_1^1=(2,1)$.

Index the partitions \eqref{eq:latticecounter} by their cyclic class:
\begin{equation}\label{eq:Pclass}
 \mathcal P_d(r)=\{K_d(y):y\in\Delta_d,\ F_d(y)\equiv r\pmod{2^d}\}.
\end{equation}
Thus $\mathcal P_d(r)=\Pi_x$ whenever $F_d(x)\equiv r$, and
$|\mathcal P_d(r)|=c_d(r)$. In the language of
Section~\ref{sec:black}, $\mathcal P_d(F_d(v))$ is the \emph{black
partition} of the unit cube $C_v$, the family of faces along which black
cubes touch $C_v$, translated to $K_d$. For a family of words, appending
a symbol means appending it to each word. Splitting the last displacement
coordinate according to its value $0$, $1$ or $-1$, whose face symbols are
$*$, $0$ and $1$, gives the face version of \eqref{eq:cyclicstern}:
\begin{align}
 \mathcal P_{d+1}(2r)&=\mathcal P_d(r)*,\label{eq:peven}\\
 \mathcal P_{d+1}(2r+1)&=\mathcal P_d(r)0\ \cup\ \mathcal P_d(r+1)1.
 \label{eq:podd}
\end{align}
Indeed, $F_{d+1}(y,s)=2F_d(y)+s$. At dimension zero we use the one-part
partition consisting of the empty word, and $c_0\equiv1$.

\subsection{Linking partitions and translated Kepler cubes}
The following lemma links the geometry of the tiling with the
compatibility condition \eqref{eq:compatible}. The parts are prescribed by
the tiling; no choice based on the numerical labels is required.

\begin{lemma}\label{lem:link}
For every $d\ge1$ and every integer $t$, the partition
$\mathcal L_d^t=\mathcal P_d(2t)$ has the following properties:
\begin{enumerate}[(1)]
\item every part has its last coordinate free;
\item if $u,v$ are antipodes in a part, then $a_d^t(u)=a_d^{t+1}(v)$;
\item its diagonal-label multisets satisfy
\begin{equation}\label{eq:spectrum}
 \Dg(a_d^t;\mathcal L_d^t)=\Ed_d(a_d^t)
 =\Ed_d(a_d^{t+1})=\Dg(a_d^{t+1};\mathcal L_d^t).
\end{equation}
\end{enumerate}
We call $\mathcal L_d^t$ a \emph{linking partition} for the pair
$(a_d^t,a_d^{t+1})$.
\end{lemma}
\begin{proof}
If $u,v$ are antipodes of $K_d(y)$, then $u+v=\ones-y$. Consequently, for
a part of $\mathcal P_d(r)$,
\begin{equation}\label{eq:antipodes}
 F_d(u)+F_d(v)\equiv-1-r\pmod{2^d}.
\end{equation}
For $r=2t$, the periodicity and evenness of $c_d$ give
$a_d^{t+1}(v)=c_d(-t-F_d(u))=a_d^t(u)$, which is (2). Property (1) is
\eqref{eq:peven}: every part of $\mathcal P_d(2t)$ has the form $S*$.

Write a diagonal of $S*$ uniquely as $\{u0,v1\}$, with $v=\tau_Su$.
Equation \eqref{eq:antipodes}, together with $F_d(u1)=F_d(u0)+1$, gives
\[
 a_d^t(u1)=a_d^t(v1),\qquad a_d^{t+1}(u0)=a_d^{t+1}(v0).
\]
Thus, in the labeling $a_d^t$, the diagonal has the same endpoint labels as
the edge $\{u0,u1\}$, and in $a_d^{t+1}$ the same as the edge
$\{v0,v1\}$. These assignments are bijections between diagonals and
$d$-edges, since $u$ runs through every vertex of every part $S$ and the
antipodal maps permute those vertices. Finally, property (2) applied to
the antipodes $u0,v1$ and $v1,u0$ shows that each diagonal carries the
same unordered pair in both labelings. This proves \eqref{eq:spectrum}.
\end{proof}

\begin{remark}\label{rem:three}
The odd classes relate three consecutive labelings. On every part of
$\mathcal P_d(2t-1)$, antipodes satisfy
\[
 a_d^t(u)=a_d^t(v),\qquad a_d^{t-1}(u)=a_d^{t+1}(v);
\]
for singleton parts the antipodal map is the identity. This is not needed
below.
\end{remark}

\begin{theorem}\label{thm:allcubes}
Every labeling $a_d^t$ is a Kepler cube. For $d\ge2$, top generating
partitions are
\begin{equation}\label{eq:top}
 \mathcal G_d^{2n}=\mathcal P_{d-1}(2n)0,\qquad
 \mathcal G_d^{2n+1}=\mathcal P_{d-1}(2n)1.
\end{equation}
In particular, for every $x\in\Z^d$, the counter labeling
$u\mapsto\ell_d(x+u)$ and the containment labeling $u\mapsto m_d(x+u)$
of $K_d$ are Kepler cubes.
\end{theorem}
\begin{proof}
We induct on $d$; the base case is the edge labeled $1,2$. For the
induction step put $A=a_{d-1}^n$, $B=a_{d-1}^{n+1}$ and
$\mathcal L=\mathcal P_{d-1}(2n)$. If $u,v$ are antipodes in a part of
$\mathcal L$, write $p=A(u)$ and $q=A(v)$. Lemma~\ref{lem:link} gives
$B(v)=p$ and $B(u)=q$. By \eqref{eq:cubrec}, in the even cube $(A,A+B)$ the
labels at $u0,v0$ are $p,q$, and those at $u1,v1$ are both $p+q$. In the
odd cube $(A+B,B)$ the same holds with generating facet $1$ and with $p,q$
interchanged. This proves \eqref{eq:generator} for the partitions
\eqref{eq:top}, whose parts are positive-dimensional by
Lemma~\ref{lem:link}(1).

Equation~\eqref{eq:spectrum} identifies the diagonal-label multiset of
$\mathcal L$ with the last-coordinate edge-label multiset of $A$ and also
of $B$. This is exactly \eqref{eq:compatible} on the generating facet, in
either parity. On the generating facet the labeling is $A$ in the even case
and $B$ in the odd case, so the induction hypothesis supplies all smaller
faces of the flag and their generating partitions; adjoining $H_d=K_d$
completes the flag. The last assertion follows from
\eqref{eq:translations}.
\end{proof}

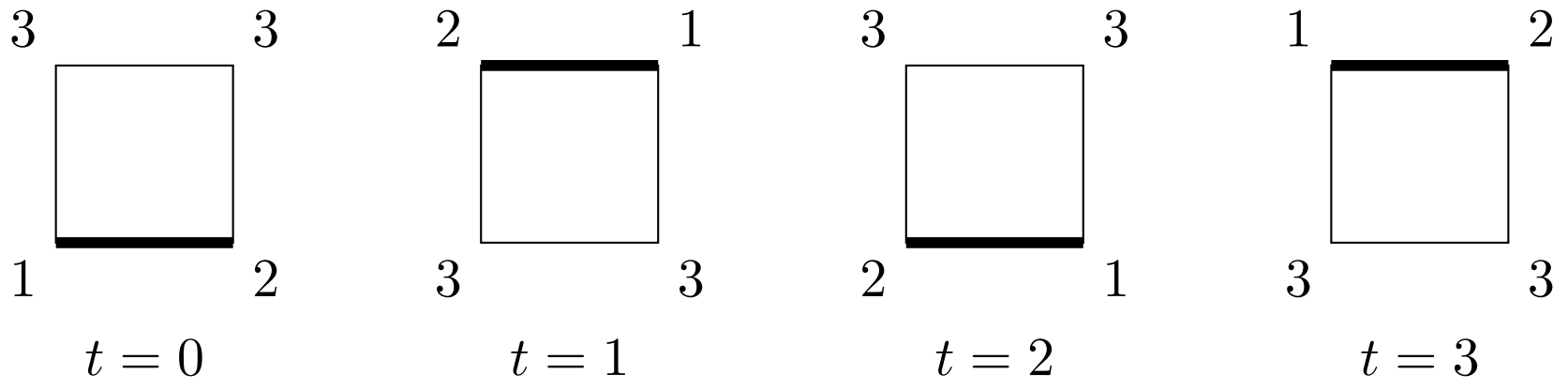
\begin{figure}[ht]
\centering
\begin{tikzpicture}[scale=1.25]
\foreach \ph/\a/\b/\c/\d/\g in {0/1/2/3/3/0,1/3/3/2/1/1,2/2/1/3/3/0,3/3/3/1/2/1}
{\begin{scope}[xshift=\ph*2.4cm]
 \draw (0,0) rectangle (1,1);
 \draw[line width=2.2pt] (0,\g)--(1,\g);
 \node[below left] at (0,0) {$\a$};
 \node[below right] at (1,0) {$\b$};
 \node[above left] at (0,1) {$\c$};
 \node[above right] at (1,1) {$\d$};
 \node at (0.5,-0.65) {$t=\ph$};
\end{scope}}
\end{tikzpicture}
\caption{The four labelings $a_2^t$ of $K_2$ (first coordinate
horizontal). The bold edge is the top generating partition
$\mathcal G_2^t$: its labels $p,q$ add up to the labels $p+q=3$ of the
opposite edge, and the two vertical edges carry the first Kepler level
$1/3$, $2/3$. Passing from $t$ to $t+1$ rotates the ring $c_2$ of
Figure~\ref{fig:rings}.}
\label{fig:squares}
\end{figure}

The occurrence of the Kepler tree in the translated labelings is now a
consequence of Proposition~\ref{prop:treeconsequence}, rather than a
property that must be verified separately. Figure~\ref{fig:squares} shows
the case $d=2$.

\begin{corollary}\label{cor:count}
Each of the families $\{u\mapsto\ell_d(x+u)\}$ and $\{u\mapsto m_d(x+u)\}$,
$x\in\Z^d$, consists of exactly $2^d$ distinct Kepler cubes, each
represented once by $x\in K_d$.
\end{corollary}
\begin{proof}
Equation~\eqref{eq:cyclicstern} shows inductively that $c_d$ equals $1$
only at residue zero. Thus the unique vertex labeled $1$ in $a_d^t$ has
binary value $-t$, so distinct phases give distinct labelings. Now use
\eqref{eq:translations} and the bijection $F_d:K_d\to\Z/2^d\Z$.
\end{proof}

\begin{example}\label{ex:product}
The Kepler property does not hold for every cube tiling lattice. For
$\Lambda=2\Z^2$, Proposition~\ref{prop:product} gives
\[
\begin{array}{c|rrrr}
 u&00&01&10&11\\\hline
 \ell_\Lambda(u)&1&2&2&4
\end{array}
\]
The two edges in direction $2$ carry $\{1,2\}$ and $\{2,4\}$, both of
ratio $1/2$, instead of the first Kepler level $1/3,2/3$. By
Proposition~\ref{prop:treeconsequence}, this labeling is not a Kepler cube.
Directly: the only possible $H_1$ is ${*}0$, with labels $1,2$, whereas the
opposite edge has labels $2,4$ and the sum rule would require $3,3$.
\end{example}

\subsection{Unique laminations and Calkin--Wilf levels}
A face partition of $K_d$ is a \emph{$k$-lamination} if every word has a
fixed symbol $0$ or $1$ in coordinate $k$. Geometrically, all parts lie in
the two facets orthogonal to $e_k$, so the cube is sandwiched in direction
$k$ between the black cubes touching it. The direction is \emph{unique} if
every other coordinate is free in at least one part. Let
$\mathcal P^k_i$ denote the subfamily of parts with $k$th symbol $i$,
i.e.\ the parts lying in the facet $x_k=i$.

\begin{theorem}\label{thm:lamination}
Every black partition of a white cube has a unique lamination direction.
More precisely, if $0<t<2^d$, then $\mathcal P_d(t)$ is laminated in the
unique direction $k=d-\nu_2(t)$. For each fixed $1\le k\le d$, exactly
$2^{k-1}$ of the partitions $\mathcal P_d(t)$, $0<t<2^d$, are
$k$-laminations, and their ratios $|\mathcal P^k_1|/|\mathcal P^k_0|$ are
precisely the fractions of level $k-1$ of the Calkin--Wilf tree, each
occurring once.
\end{theorem}
\begin{proof}
For $0<t<2^d$ write $t=2^{d-k}(2r+1)$, with $0\le r<2^{k-1}$.
Equations~\eqref{eq:peven}--\eqref{eq:podd} give
\begin{equation}\label{eq:laminationform}
 \mathcal P_d(t)=
 \bigl(\mathcal P_{k-1}(r)0\ \cup\ \mathcal P_{k-1}(r+1)1\bigr)*^{d-k}.
\end{equation}
So coordinate $k$ is fixed in every part and the last $d-k$ coordinates
are free. We prove by induction on $k$ that each of the coordinates
$1,\ldots,k-1$ is free in some part. For $k=1$ there is nothing to prove.
For $k>1$, one of the indices $r,r+1$ is odd modulo $2^{k-1}$ and the other
even. By induction, the odd-indexed partition is laminated in the unique
direction $k-1$, so each of the coordinates $1,\ldots,k-2$ is free in some
of its parts; by \eqref{eq:peven}, the even-indexed partition has
coordinate $k-1$ free in every part. Hence no coordinate other than $k$ is
fixed in every part of \eqref{eq:laminationform}.

The ratio in \eqref{eq:laminationform} is $c_{k-1}(r+1)/c_{k-1}(r)$. By
\eqref{eq:cyclicstern}, the cyclically adjacent pairs
$(c_j(r),c_j(r+1))$, $0\le r<2^j$, are replaced in the next dimension by
$(p,p+q)$ and $(p+q,q)$ for each pair $(p,q)$. Starting from $c_0\equiv1$,
they are therefore the Calkin--Wilf levels as ordered pairs.
Reversing the coordinates permutes each level (see the proof of
Proposition~\ref{prop:mixed}), so the reciprocal ratios give the same
level, each fraction once. Distinct residues give distinct partitions,
since $y\mapsto K_d(y)$ is injective.
\end{proof}

In Figure~\ref{fig:black} ($d=2$), a white square in an odd row has the
partition $\mathcal P_2(1)=\{{*}0,01,11\}$ or
$\mathcal P_2(3)=\{00,10,{*}1\}$, laminated vertically with ratios $2/1$
and $1/2$. A white square in an even row has
$\mathcal P_2(2)=\{0{*},1{*}\}$, laminated horizontally with ratio $1/1$.
These are levels $1$ and $0$ of the Calkin--Wilf tree.

\subsection{Illustrations in dimension three}
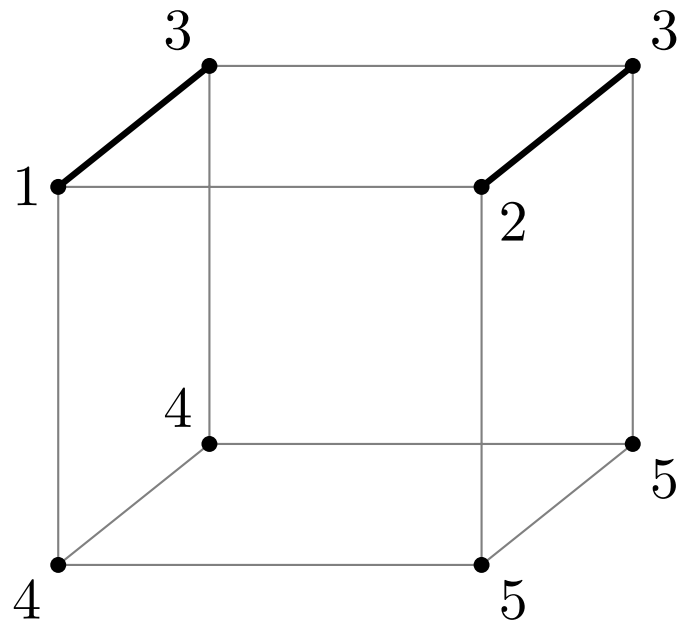
\begin{figure}[ht]
\centering
\begin{tikzpicture}[x={(2.8cm,0cm)},y={(1cm,0.8cm)},z={(0cm,2.5cm)}]
\coordinate (O) at (0,0,0);\coordinate (A) at (1,0,0);
\coordinate (B) at (1,1,0);\coordinate (C) at (0,1,0);
\coordinate (D) at (0,0,1);\coordinate (E) at (1,0,1);
\coordinate (F) at (1,1,1);\coordinate (G) at (0,1,1);
\draw[gray] (O)--(A)--(B)--(C)--cycle
 (D)--(E)--(F)--(G)--cycle (O)--(D) (A)--(E) (B)--(F) (C)--(G);
\draw[very thick] (D)--(G) (E)--(F);
\foreach \p/\v/\pos in {O/4/below left,A/5/below right,B/5/below right,C/4/above left,
D/1/left,E/2/below right,F/3/above right,G/3/above left}{
 \fill (\p) circle (1.5pt);
 \node[\pos,inner sep=3pt] at (\p) {$\v$};
}
\end{tikzpicture}
\caption{The containment labeling $m_3$ on the cube at the origin. The
first coordinate runs to the right, the second along the oblique edges,
and the third vertically. The bold edges $0{*}1$ and $1{*}1$ are the parts
of the top generating partition.}
\label{fig:cube}
\end{figure}

At the origin, the containment labeling has phase $-1$ by
\eqref{eq:translations}. In dimension three it is
\[
\begin{array}{c|rrrrrrrr}
 u&000&001&010&011&100&101&110&111\\\hline
 m_3(u)=a_3^{-1}(u)&4&1&4&3&5&2&5&3
\end{array}
\]
The top generating partition is $\mathcal P_2(2)1=\{0{*}1,1{*}1\}$. Its
two diagonals carry the pairs $\{1,3\}$ and $\{2,3\}$, while the opposite
labels are $4,4$ and $5,5$. The vertical edges therefore give
$1/4,3/4,2/5,3/5$, the second Kepler level. A compatible flag is
$H_1={*}01$, $H_2={*}{*}1$, $H_3={*}{*}{*}$; see Figure~\ref{fig:cube}.

The \emph{diagonal staircase} of unit cubes $C_{k\ones}$,
$0\le k\le2^d$, runs through all black partitions, in the order
\[
 \mathcal P_d(0),\ \mathcal P_d(-1),\ \ldots,\ \mathcal P_d(-2^d),
\]
since $F_d(k\ones)\equiv-k$. Only the first and the last cube are black;
their partition is $\{{*}^d\}$, and every other partition occurs once. The
point $k\ones$ is a vertex of a unique black cube, and its local
coordinate there has binary value $-k\bmod2^d$. Reading these coordinates
for $k=2^d,2^d-1,\ldots,1$ lists all binary words of length $d$ in
increasing order.

%======================================================================
\section{Reflections and cyclic profiles}\label{sec:profiles}
%======================================================================
\subsection{Coordinate reflections}
Let $D_\delta=\diag(\delta_1,\ldots,\delta_d)$ with $\delta_i\in\{-1,1\}$.
Coordinate reflections preserve $\Delta_d$, so for every cube tiling
lattice $\Lambda$
\begin{equation}\label{eq:reflection}
 \ell_{D_\delta\Lambda}(x)=\ell_\Lambda(D_\delta x).
\end{equation}
The reflected lattice tiles, because $D_\delta K_d$ is a translate of
$K_d$. If $\Lambda=A\Z^d$ is standard, then $D_\delta\Lambda=D_\delta
AD_\delta\Z^d$, and $D_\delta AD_\delta$ is upper triangular with diagonal
$2$; normalizing its columns shows that $D_\delta\Lambda$ is standard.
Cocyclicity is preserved by the isomorphism $D_\delta$.

\begin{corollary}\label{cor:profile}
For every $i$ and every integer $k$,
$\ell_{D_\delta\Lambda}(ke_i)=\ell_\Lambda(ke_i)$.
\end{corollary}
\begin{proof}
Combine \eqref{eq:reflection} with the evenness of $\ell_\Lambda$.
\end{proof}

In terms of normalized weights, $D_\delta\Lambda_w$ is the kernel of
$x\mapsto\sum_i\delta_iw_ix_i$; normalizing the last weight to $1$, its
normalized weights are
\begin{equation}\label{eq:reflectedweights}
 w_i'\equiv\delta_d\delta_iw_i\pmod N.
\end{equation}

\begin{theorem}\label{thm:orbit}
For $d\ge2$, every orbit of a standard cocyclic cube tiling lattice under
coordinate reflections has exactly $2^{d-2}$ members. The number of
orbits is $2^{(d-2)(d-3)/2}$.
\end{theorem}
\begin{proof}
By \eqref{eq:reflectedweights} and the uniqueness in
Theorem~\ref{thm:weights}, $D_\delta\Lambda=\Lambda$ if and only if
$\delta_d\delta_iw_i\equiv w_i$ for all $i$. By
Corollary~\ref{cor:cyclicbinary}, $\ord(w_i)=2^i$. For $i\ge2$ we have
$w_i\ne-w_i$, so $\delta_i=\delta_d$; for $i=1$ both signs are allowed,
as $w_1=-w_1$. There are exactly four stabilizing sign vectors, so every
orbit has $2^d/4$ elements. Dividing the count of
Corollary~\ref{cor:cyclicbinary} by $2^{d-2}$ gives the number of orbits.
\end{proof}

\begin{remark}\label{rem:sternorbit}
By Remark~\ref{rem:signedweights}, the signed lattices
$\Lambda_d^\varepsilon$ are exactly the standard lattices with
$w_i\equiv\pm2^{d-i}$. By \eqref{eq:reflectedweights}, they form the
reflection orbit of $\Lambda_d^b$, whose weights are $2^{d-i}$. This
explains Corollary~\ref{cor:mixedlattices}: the $2^{d-2}$ lattices
\eqref{eq:gradual} form one orbit.
\end{remark}

\subsection{The cyclic profile and its Fourier transform}
For a standard cocyclic lattice $\Lambda$ with normalized weights $w$ and
$N=2^d$, define its \emph{(ordered) cyclic profile}
\[
 C_\Lambda=(c_k)_{k\in\ZN},\qquad c_k=\ell_\Lambda(ke_d).
\]
The points $ke_d$, $k\in\ZN$, represent all cosets. Since
$\ell_\Lambda(x)=|\{y\in\Delta_d:x+y\in\Lambda\}|$ and $\Delta_d=-\Delta_d$,
\begin{equation}\label{eq:profilecount}
 \ell_\Lambda(x)=c_{w\cdot x},\qquad
 c_k=\Bigl|\Bigl\{y\in\Delta_d:\sum_iw_iy_i\equiv k\pmod N\Bigr\}\Bigr|.
\end{equation}
Thus the profile determines the whole counter. Each $y\in\Delta_d$
contributes a monomial, and
\begin{equation}\label{eq:polynomial}
 \sum_{k=0}^{N-1}c_kz^k\equiv\prod_{i=1}^d\bigl(1+z^{w_i}+z^{-w_i}\bigr)
 \pmod{z^N-1}.
\end{equation}
For the Stern orbit, $w_i=\pm2^{d-i}$, and since $z^{1-2^d}\equiv z$,
\[
 \sum_{k=0}^{N-1}c_kz^k\equiv z\prod_{j=0}^{d-1}
 \bigl(1+z^{2^j}+z^{2^{j+1}}\bigr)\pmod{z^N-1},
\]
a truncation of the classical generating function
$\prod_{j\ge0}(1+z^{2^j}+z^{2^{j+1}})=\sum_{n\ge0}s(n+1)z^n$ of Stern's
sequence; see \cite{Northshield}.

We evaluate \eqref{eq:polynomial} at the $N$th roots of unity. For
$a\in\ZN$ let $\chi_a(x)=\exp(2\pi i\,ax/N)$, $x\in\ZN$. We work in the
algebra $\C^{\ZN}$ of complex functions on $\ZN$ with pointwise
multiplication. The characters form an orthonormal basis for the inner
product $\langle f,g\rangle=N^{-1}\sum_xf(x)\overline{g(x)}$, and
$\chi_a\chi_b=\chi_{a+b}$. The Fourier transform of the profile is
\begin{equation}\label{eq:fourier}
 \widehat C_\Lambda(x)=\sum_{k\in\ZN}c_k\chi_k(x)
 =\sum_{y\in\Delta_d}\chi_{w\cdot y}(x)=\Phi_w(x),\qquad
 \Phi_w=\prod_{i=1}^d\bigl(1+\chi_{w_i}+\chi_{-w_i}\bigr).
\end{equation}
Since the Fourier transform is invertible, $C_\Lambda$ and $\Phi_w$
determine each other.

\subsection{Profiles determine lattices up to reflections}
\begin{lemma}\label{lem:unit}
Let $N$ be a power of two. For every $a\in\ZN$, the function
$1+\chi_a+\chi_{-a}$ has no zeros on $\ZN$; hence it is invertible in
$\C^{\ZN}$.
\end{lemma}
\begin{proof}
If $1+t+t^{-1}=0$ for $t=\chi_a(x)$, then $t^2+t+1=0$, so $t$ is a
primitive third root of unity. But $t$ is an $N$th root of unity and $N$ is
a power of two.
\end{proof}

\begin{theorem}\label{thm:unique}
Let $N=2^d$, and let $w_i,v_i\in\ZN$ have order exactly $2^i$ for
$1\le i\le d$. Then $\Phi_w=\Phi_v$ if and only if $v_i=\pm w_i$ for all
$i$, with independent signs.
\end{theorem}
\begin{proof}
Changing the sign of a weight does not change its factor, which proves one
implication. We prove the converse by induction on $d$; the case $d=1$ is
immediate since then $w_1=v_1=1$. Let $d\ge2$. Since $\ord(w_i)=2^i<N$
for $i<d$, the elements $w_1,\ldots,w_{d-1}$ are even, whereas $w_d$ is
odd; the same holds for $v$. Put
\[
 A_w=\prod_{i<d}\bigl(1+\chi_{w_i}+\chi_{-w_i}\bigr),\qquad
 A_v=\prod_{i<d}\bigl(1+\chi_{v_i}+\chi_{-v_i}\bigr),
\]
and let $E$ and $O$ be the spans of the characters with even and with odd
frequencies, respectively. Then $\C^{\ZN}=E\oplus O$, orthogonally. All
frequencies in $A_w$ and $A_v$ are even, while
$(\chi_{w_d}+\chi_{-w_d})A_w$ and $(\chi_{v_d}+\chi_{-v_d})A_v$ lie in $O$.
Hence
\[
 \Phi_w=A_w+(\chi_{w_d}+\chi_{-w_d})A_w,\qquad
 \Phi_v=A_v+(\chi_{v_d}+\chi_{-v_d})A_v
\]
are the decompositions along $E\oplus O$, and $\Phi_w=\Phi_v$ implies
$A_w=A_v$.

The identification $\chi^{(N)}_{2a}\leftrightarrow\chi^{(N/2)}_a$ maps
$E$ isomorphically, as an algebra, onto the function algebra of
$\Z_{N/2}$: the functions in $E$ are exactly the $N/2$-periodic ones.
Halving the even residues $w_i,v_i$, $i<d$, preserves their orders, so
the induction hypothesis gives $v_i/2\equiv\pm w_i/2\pmod{N/2}$, that is,
$v_i=\pm w_i$ for $i<d$.

Finally, by Lemma~\ref{lem:unit}, the common function $A_w=A_v$ has no
zeros. Cancelling it in $\Phi_w=\Phi_v$ gives
$\chi_{w_d}+\chi_{-w_d}=\chi_{v_d}+\chi_{-v_d}$, and the linear
independence of characters yields $v_d=\pm w_d$.
\end{proof}

\begin{theorem}\label{thm:main}
Two standard cocyclic cube tiling lattices in $\Z^d$ have the same ordered
cyclic profile if and only if they differ by coordinate reflections.
\end{theorem}
\begin{proof}
Let $w$ and $v$ be the normalized weights of the two lattices; by
Corollary~\ref{cor:cyclicbinary}, $\ord(w_i)=\ord(v_i)=2^i$. If the
profiles agree, then $\Phi_w=\Phi_v$ by \eqref{eq:fourier}, and
Theorem~\ref{thm:unique} gives $\delta_i\in\{-1,1\}$ with
$v_i\equiv\delta_iw_i$. Since $v_d=w_d=1$ we may take $\delta_d=1$, and
\eqref{eq:reflectedweights} shows that $\Lambda_v=D_\delta\Lambda_w$. The
converse is Corollary~\ref{cor:profile}.
\end{proof}

\begin{corollary}\label{cor:numberprofiles}
For $d\ge2$, the standard cocyclic cube tiling lattices in $\Z^d$ have
exactly $2^{(d-2)(d-3)/2}$ distinct cyclic profiles. A standard lattice
has the Stern profile \eqref{eq:sternrow} if and only if it is one of the
signed lattices $\Lambda_d^\varepsilon$.
\end{corollary}
\begin{proof}
Combine Theorems~\ref{thm:orbit} and~\ref{thm:main} with
Remark~\ref{rem:sternorbit}.
\end{proof}

\subsection{Examples}
For $d=4$ there are eight standard cocyclic lattices in two orbits. With
normalized weights $(8,4,2,1)$ and $(8,4,6,1)$, representatives are
\[
 B_4=\begin{pmatrix}2&1&1&1\\0&2&1&1\\0&0&2&1\\0&0&0&2\end{pmatrix},
 \qquad
 \begin{pmatrix}2&1&0&1\\0&2&1&0\\0&0&2&1\\0&0&0&2\end{pmatrix},
\]
with the profiles
\[
\begin{array}{l}
 1,5,4,7,3,8,5,7,2,7,5,8,3,7,4,5\quad\text{(Stern)},\\
 1,6,5,8,3,7,4,6,2,6,4,7,3,8,5,6.
\end{array}
\]
Both sum to $3^4=81$, in accordance with Theorem~\ref{thm:sum}. The
following table compares the counts of Theorem~\ref{thm:orbit} with an
exhaustive computation based on \eqref{eq:polynomial}; the last column is
discussed in Remark~\ref{rem:multisets}.
\begin{center}
\begin{tabular}{rrrr}
\toprule
$d$&lattices $2^{(d-1)(d-2)/2}$&profiles $=$ orbits&distinct multisets\\
\midrule
2&1&1&1\\
3&2&1&1\\
4&8&2&2\\
5&64&8&8\\
6&1\,024&64&63\\
7&32\,768&1\,024&1\,022\\
\bottomrule
\end{tabular}
\end{center}

%======================================================================
\section{Lattices tiling by discrete boxes}\label{sec:boxes}
%======================================================================
\subsection{Partition counters for boxes}
Let $n_1,\ldots,n_d\ge2$, and let $K=K(n)$ and $\Delta=\Delta(n)$ be as in
Section~\ref{sec:bases}. For $y\in\Delta$,
\[
 K(y)=K\cap(K-y)=\prod_{i=1}^d
 \{\max(0,-y_i),\ldots,\min(n_i-1,n_i-1-y_i)\}.
\]
The map $y\mapsto K(y)$ is injective, and its image consists of the
subboxes each of whose coordinate intervals contains $0$ or $n_i-1$.
Suppose that $K+T$ partitions $\Z^d$, and define
$\ell_T(x)=|T\cap(x+\Delta)|$. Since $-K=K-(n-\ones)$, the proof of
Theorem~\ref{thm:sum} gives
\[
 \sum_{a\in K+z}\ell_T(a)=\prod_{i=1}^d(2n_i-1)\qquad(z\in\Z^d),
\]
and for a lattice the displacement formula of
Proposition~\ref{prop:displacement} holds with the representative map of
$K$.

From now on let $\Lambda=A\Z^d$, where $A$ is upper triangular with
diagonal $n_1,\ldots,n_d$ having common prime support, and assume that
$\Lambda$ is cocyclic. If the coordinates are reordered to make a lattice
triangular, the side lengths must be reordered with them: the diagonal
entry in coordinate $j$ is the side length of the box in that coordinate.
Let $w$ be the normalized weights of Theorem~\ref{thm:weights}, so that
$\ord(w_i)=N_i$, and define the ordered cyclic profile
$C_\Lambda=(c_k)_{k\in\ZN}$, $c_k=\ell_\Lambda(ke_d)$, as before. Exactly
as in \eqref{eq:profilecount}--\eqref{eq:fourier},
\begin{equation}\label{eq:boxfourier}
 \widehat C_\Lambda=\prod_{i=1}^d\theta_{n_i,w_i},\qquad
 \theta_{n,a}=\sum_{r=-(n-1)}^{n-1}\chi_{ra}.
\end{equation}
The factors $\theta_{n,a}$ are symmetric Dirichlet kernels;
$\theta_{2,a}=1+\chi_a+\chi_{-a}$.

\begin{lemma}\label{lem:invertible}
Let $n\ge2$ be divisible by every prime divisor of $N$. Then, for every
$a\in\ZN$, the function $\theta_{n,a}$ has no zeros on $\ZN$.
\end{lemma}
\begin{proof}
Let $t=\chi_a(x)$. If $t=1$, the value is $2n-1\ne0$. Otherwise
\[
 \theta_{n,a}(x)=t^{-(n-1)}\,\frac{1-t^{2n-1}}{1-t},
\]
so a zero would force $\ord(t)$ to divide both $N$ and $2n-1$. Every prime
divisor of $N$ divides $n$ and hence not $2n-1$, so $\gcd(N,2n-1)=1$ and
$t=1$, a contradiction.
\end{proof}

For $n_i=2$ this is Lemma~\ref{lem:unit}: the number $3=2n_i-1$ is the
only obstruction, and it is excluded because $N$ is a power of two. Under
the common prime support no separate assumption such as $3\nmid N$ is
needed.

\subsection{Separation of the factors}
Every $w_i$ with $\ord(w_i)=N_i$ can be written uniquely as
\[
 w_i\equiv M_ia_i\pmod N,\qquad a_i\in(\Z/N_i\Z)^\times.
\]
For $1\le i\le d$, define the \emph{rigidity group}
\[
 \mathcal R_i=\bigl\{r\in(\Z/N_i\Z)^\times:\ \theta_{n_i,r}=\theta_{n_i,1}
 \text{ as functions on }\Z/N_i\Z\bigr\},
\]
where the characters in $\theta_{n_i,r}$ are those of $\Z/N_i\Z$. It is a
subgroup containing $\pm1$, because $\theta_{n,ra}(x)=\theta_{n,a}(rx)$.

\begin{theorem}\label{thm:separation}
Let $w$ and $v$ satisfy $\ord(w_i)=\ord(v_i)=N_i$ for all $i$, and write
$w_i\equiv M_ia_i$ and $v_i\equiv M_ib_i$ as above. Then
\[
 \prod_{i=1}^d\theta_{n_i,w_i}=\prod_{i=1}^d\theta_{n_i,v_i}
 \quad\Longleftrightarrow\quad
 b_ia_i^{-1}\in\mathcal R_i\quad(1\le i\le d).
\]
\end{theorem}
\begin{proof}
For $i<d$, $\ord(w_i)$ divides $N/n_d$, so $w_i$ and $v_i$ are multiples of
$n_d$ in $\ZN$; hence all frequencies of
$A_w=\prod_{i<d}\theta_{n_i,w_i}$ and of the analogous $A_v$ are multiples
of $n_d$. On the other hand $w_d=M_da_d=a_d$ is a unit modulo $N$, and
hence modulo $n_d$. In $\theta_{n_d,w_d}=\sum_{|r|<n_d}\chi_{rw_d}$, the
frequency $rw_d$ is divisible by $n_d$ only for $r=0$; the same holds for
$v_d$. Orthogonal projection onto the span of the characters with
frequencies divisible by $n_d$ therefore gives $A_w=A_v$. By
Lemma~\ref{lem:invertible} this function has no zeros, and cancelling it
gives $\theta_{n_d,w_d}=\theta_{n_d,v_d}$.

The span of the characters with frequencies divisible by $n_d$ is the
algebra of $(N/n_d)$-periodic functions, identified with the function
algebra of $\Z/(N/n_d)\Z$ by dividing frequencies by $n_d$. Under this
identification $A_w=A_v$ is the same problem for $n_1,\ldots,n_{d-1}$,
with weights $w_i/n_d\equiv(M_i/n_d)a_i$. Iterating, we obtain for every $i$
the equality $\theta_{n_i,a_i}=\theta_{n_i,b_i}$ of functions on
$\Z/N_i\Z$. The substitution $x\mapsto a_i^{-1}x$ turns it into
$\theta_{n_i,b_ia_i^{-1}}=\theta_{n_i,1}$, that is,
$b_ia_i^{-1}\in\mathcal R_i$.

Conversely, $\theta_{n_i,w_i}$ is the pullback to $\ZN$ of the function
$\theta_{n_i,a_i}$ on $\Z/N_i\Z$, since $\chi_{M_ia}^{(N)}(x)=
\chi_a^{(N_i)}(x\bmod N_i)$. Equality of the local factors therefore
implies equality of the products.
\end{proof}

\subsection{The rigidity groups}
The groups $\mathcal R_i$ can be determined completely.

\begin{proposition}\label{prop:rigidity}
Let $n_1,\ldots,n_d\ge2$ be arbitrary. Then:
\begin{enumerate}[(a)]
\item $\mathcal R_1=(\Z/n_1\Z)^\times$;
\item if $i\ge2$ and $N_{i-1}=2$, then $\mathcal R_i=(\Z/2n_i\Z)^\times$;
\item if $i\ge2$ and $N_{i-1}\ge3$, then $\mathcal R_i=\{\pm1\}$.
\end{enumerate}
\end{proposition}
\begin{proof}
(a) On $\Z/n_1\Z$, the integers $-(n_1-1),\ldots,n_1-1$ cover the residue
$0$ once and every other residue twice. Hence $\theta_{n_1,1}(0)=2n_1-1$
and $\theta_{n_1,1}(x)=2\sum_{s}\chi_1(sx)-1=-1$ for $x\ne0$. Since
$\theta_{n_1,r}(x)=\theta_{n_1,1}(rx)$ and $rx=0$ only for $x=0$, every unit
$r$ belongs to $\mathcal R_1$.

For $i\ge2$ write $n=n_i$ and $m=N_{i-1}\ge2$, so that $N_i=mn\ge2n$.
The integers in $S=\{-(n-1),\ldots,n-1\}$ are then distinct modulo $mn$,
and $rS$ consists of distinct residues for a unit $r$. By the linear
independence of characters, $\theta_{n,r}=\theta_{n,1}$ if and only if
$rS=S$ in $\Z/mn\Z$.

(b) If $m=2$, the complement of $S$ in $\Z/2n\Z$ is $\{n\}$. Every unit
$r$ is odd and so fixes $n$; hence $rS=S$.

(c) Let $m\ge3$ and $r\in\mathcal R_i$. Since $1\in S$, we have $r\in S$,
so $r\equiv s$ with $0<|s|\le n-1$; replacing $r$ by $-r$, we may assume
$1\le s\le n-1$. Suppose $s\ge2$; then $n\ge3$. Let $k=\lceil n/s\rceil$.
Then $2\le k\le\lceil n/2\rceil\le n-1$, so $k\in S$, and
$n\le ks\le n+s-1\le2n-2$. Modulo $mn$, the set $S$ is
$\{0,\ldots,n-1\}\cup\{mn-n+1,\ldots,mn-1\}$, and
$2n-2<mn-n+1$ because $m\ge3$. Hence $rk\equiv ks\notin S$, a
contradiction. So $s=1$ and $r=\pm1$.
\end{proof}

\begin{theorem}\label{thm:boxclass}
Let $d\ge2$ and let $n_1,\ldots,n_d\ge2$ have common prime support.
Consider the cocyclic lattices with an upper triangular basis with diagonal
$n_1,\ldots,n_d$.
\begin{enumerate}[(a)]
\item Two such lattices with normalized weights $w_i\equiv M_ia_i$ and
$v_i\equiv M_ib_i$ have the same ordered cyclic profile if and only if
$b_ia_i^{-1}\in\mathcal R_i$ for $1\le i<d$.
\item The ordered cyclic profile determines each such lattice up to
coordinate reflections if and only if $\mathcal R_i=\{\pm1\}$ for
$1\le i<d$. This holds if and only if
\[
 n_1\in\{3,4,6\},\qquad\text{or}\qquad n_1=2\ \text{and}\ (d=2\ \text{or}\
 n_2=2).
\]
\end{enumerate}
\end{theorem}
\begin{proof}
(a) The profiles agree if and only if the Fourier products
\eqref{eq:boxfourier} agree. Apply Theorem~\ref{thm:separation}, noting
that $a_d=b_d=1$.

(b) Suppose $\mathcal R_i=\{\pm1\}$ for $i<d$. If the profiles agree, then
$b_i\equiv\pm a_i\pmod{N_i}$, so $v_i\equiv\pm w_i\pmod N$ for $i<d$, and
$\Lambda_v=D_\delta\Lambda_w$ by \eqref{eq:reflectedweights}, whose proof
applies verbatim. Conversely, reflections preserve $\Delta$, and
$\ell_\Lambda$ is even, so reflected lattices have the same profile.

Now suppose that some $\mathcal R_i$, $i<d$, contains $r\ne\pm1$. Let
$w_j=M_j$ for all $j$, and let $v$ agree with $w$ except that
$v_i=M_ir\bmod N$. All orders are unchanged, so by
Proposition~\ref{prop:converse} both $w$ and $v$ are the normalized
weights of cocyclic lattices with diagonal $n_1,\ldots,n_d$. By (a) they
have the same profile. They are not related by a reflection, because
$v_i\not\equiv\pm w_i\pmod N$, i.e.\ $r\not\equiv\pm1\pmod{N_i}$.

Finally, we use Proposition~\ref{prop:rigidity}. The group
$(\Z/n_1\Z)^\times$ equals $\{\pm1\}$ if and only if $\varphi(n_1)\le2$,
that is, $n_1\in\{2,3,4,6\}$. If $n_1\ge3$, then $N_{i-1}\ge3$ for all
$i\ge2$, and all $\mathcal R_i$ with $i\ge2$ are $\{\pm1\}$. If $n_1=2$ and
$d\ge3$, then $\mathcal R_2=(\Z/2n_2\Z)^\times$, which is $\{\pm1\}$ if
and only if $\varphi(2n_2)=2$, i.e.\ $n_2=2$ (the value $n_2=3$ is excluded
by the common prime support); the remaining $\mathcal R_i$, $i\ge3$, are
$\{\pm1\}$ since $N_{i-1}\ge4$.
\end{proof}

\begin{corollary}\label{cor:boxprofiles}
Under the hypotheses of Theorem~\ref{thm:boxclass}, the number of distinct
ordered cyclic profiles is
\[
 \prod_{i=1}^{d-1}\frac{\varphi(N_i)}{|\mathcal R_i|},
\]
where $|\mathcal R_1|=\varphi(n_1)$, $|\mathcal R_2|=\varphi(2n_2)$ if
$n_1=2$, and $|\mathcal R_i|=2$ in all other cases with $i\ge2$.
\end{corollary}
\begin{proof}
By Proposition~\ref{prop:converse}, the lattices correspond to the tuples
\[
 (a_1,\ldots,a_{d-1})\in\prod_{i<d}(\Z/N_i\Z)^\times,
\]
and by Theorem~\ref{thm:boxclass}(a) the profiles correspond to the cosets
of $\prod_{i<d}\mathcal R_i$.
\end{proof}

For $n_i=2$ this gives $\prod_{i=2}^{d-1}2^{i-2}=2^{(d-2)(d-3)/2}$, in
accordance with Corollary~\ref{cor:numberprofiles}.

\begin{remark}\label{rem:boxorbits}
By \eqref{eq:reflectedweights}, a reflection $D_\delta$ fixes the lattice
if and only if $\delta_i=\delta_d$ for every $i<d$ with $2w_i\not\equiv0$.
Since $2w_i\equiv0$ only when $N_i\le2$, that is, $i=1$ and $n_1=2$, the
reflection orbits have $2^{d-1}$ elements if $n_1\ge3$ and $2^{d-2}$
elements if $n_1=2$ ($d\ge2$).
\end{remark}

\begin{example}\label{ex:five}
Let $n_1=n_2=5$, so $N=25$ and $\mathcal R_1=(\Z/5\Z)^\times$. The four
cocyclic lattices have normalized weights $(5,1)$, $(10,1)$, $(15,1)$ and
$(20,1)$; they form two reflection orbits, $\{(5,1),(20,1)\}$ and
$\{(10,1),(15,1)\}$, but all four have the profile
\[
 1,3,3,3,3,\ 2,4,4,4,4,\ 2,4,4,4,4,\ 2,4,4,4,4,\ 2,3,3,3,3.
\]
For instance, $\ker(5x_1+x_2\bmod25)$ and $\ker(10x_1+x_2\bmod25)$ have
the Hermite normal forms
\[
 \begin{pmatrix}5&4\\0&5\end{pmatrix},\qquad
 \begin{pmatrix}5&2\\0&5\end{pmatrix},
\]
as $4\cdot5+5\equiv0$ and $2\cdot10+5\equiv0\pmod{25}$.
\end{example}

\begin{example}\label{ex:242}
Even for $2$-power sides the reconstruction may fail. Let
$(n_1,n_2,n_3)=(2,4,2)$, so $N=16$ and $\mathcal R_2=(\Z/8\Z)^\times$.
The lattices
\[
 \ker(8x_1+2x_2+x_3\bmod16)=\begin{pmatrix}2&1&1\\0&4&3\\0&0&2
 \end{pmatrix}\Z^3,\qquad
 \ker(8x_1+6x_2+x_3\bmod16)=\begin{pmatrix}2&1&1\\0&4&1\\0&0&2
 \end{pmatrix}\Z^3
\]
both have the profile $1,4,3,6,3,6,3,5,2,5,3,6,3,6,3,4$, but
$6\not\equiv\pm2\pmod{16}$, so they are not related by a reflection.
\end{example}

%======================================================================
\section{Further remarks and questions}\label{sec:questions}
%======================================================================
\begin{question}\label{q:nonlattice}
Which partitions of $\mathscr K_d$ into partitions of $K_d$ arise as
$\{\Pi_a^T:a\in K_d+z\}$ from a discrete cube tiling? How does this
collection depend on $z$ when $T$ is not a lattice?
\end{question}

For a lattice, changing $z$ merely permutes the coset representatives and
leaves the unindexed collection unchanged. No such periodicity argument is
available for a general translation set.

\begin{remark}\label{rem:multisets}
Theorem~\ref{thm:main} concerns the ordered profile, with respect to the
generator $e_d$. The multiset of its values, equivalently the multiset of
counter values on $K_d$, is a weaker invariant. An exhaustive computation
shows that it still separates the reflection orbits for $d\le5$, but not
for $d=6$ and $d=7$ (see the table in Section~\ref{sec:profiles}). For
instance, the lattices with normalized weights $(32,16,8,12,26,1)$ and
$(32,16,8,20,6,1)$ in $\Z^6$ lie in different orbits but have the same
multiset of counter values. It would be interesting to know which
invariants of a cocyclic lattice are encoded by this multiset.
\end{remark}

%======================================================================
\appendix
\section{Group transversals}\label{app:groups}
%======================================================================
The displacement construction of Section~\ref{sec:partitions} has a
general group-theoretic form. Let $G$ be a group, let
$H\trianglelefteq G$, and let $K$ contain exactly one representative of
every coset of $H$. Write $a\circ b$ for the representative of $abH$ in
$K$; an inverse in $(K,\circ)$ is denoted by $a'$, whereas $a^{-1}$ is the
inverse in $G$. Put
\[
 L=KK^{-1},\qquad K(x)=K\cap x^{-1}K\quad(x\in L),
 \qquad L_a=L\cap aH.
\]
Every $K(x)$ is nonempty.

\begin{proposition}\label{prop:group}
For each $a\in K$,
\[
 L_a=\{(a\circ u)u^{-1}:u\in K\}.
\]
The family $\Pi_a=\{K(x):x\in L_a\}$ is a partition of $K$, and
$a\circ u=xu$ for $u\in K(x)$, $x\in L_a$. Distinct elements of $L_a$
give distinct parts. If $K$ is finite, then $|\Pi_a|=|L_a|$ and
$|\Pi_{a'}|=|\Pi_a|$.
\end{proposition}
\begin{proof}
If $x=vu^{-1}\in L_a$, then $vH=auH$, so $v=a\circ u$. This proves the
formula for $L_a$, and also the translation formula, by the uniqueness of
representatives. Every $u$ lies in the part indexed by $(a\circ u)u^{-1}$.
If $u$ lies in both $K(x)$ and $K(y)$ for $x,y\in L_a$, then
$xu=a\circ u=yu$, so $x=y$. Finally, inversion maps $L_a$ bijectively onto
$L_{a'}$.
\end{proof}

The translates $hK$, $h\in H$, partition $G$. If $gH=aH$ and
$P_g=\{gK\cap hK:h\in H\}\setminus\{\varnothing\}$, then $P_g=g\Pi_a$:
indeed $gK\cap hK=gK(h^{-1}g)$, and $h^{-1}g$ ranges through $L_a$
precisely when the intersection is nonempty. For $x\in L_a$ we also have
$a\circ K(x)=xK(x)=K(x^{-1})$, so the permutation $u\mapsto a\circ u$
carries $\Pi_a$ onto $\Pi_{a'}$.

\begin{remark}\label{rem:injectivity}
The sets $L_a$, $a\in K$, always partition $L$, but the families $\Pi_a$
need not be pairwise disjoint as families of subsets. They are if
$x\mapsto K(x)$ is injective on $L$, as it is for discrete cubes and boxes;
in general one can use the labeled parts $(x,K(x))$ instead.
\end{remark}

\begin{example}
In additive notation, take $G=\Z$, $H=4\Z$ and $K=\{0,1,2,7\}$. Then
$K(5)=K(-2)=\{2\}$, while $5$ and $-2$ represent different cosets modulo
$4\Z$. Explicitly,
\[
 \Pi_1=\{\{0,1\},\{2\},\{7\}\},\qquad
 \Pi_2=\{\{0\},\{1\},\{2\},\{7\}\},
\]
which are different partitions with common members.
\end{example}

%======================================================================
\section{The direct-limit construction}\label{app:infinite}
%======================================================================
Let $\Z^{(\mathbb N)}$ be the group of finitely supported integer
sequences, and let $K_\infty$ consist of its $0$--$1$ sequences. For fixed
signs as in \eqref{eq:signed}, let
\[
 \Lambda_\infty=\langle2e_1,\ 2e_j+\varepsilon_{j-1}e_{j-1}\ (j\ge2)
 \rangle_{\Z}.
\]
Every vector can be reduced within a finite-dimensional coordinate
subspace by \eqref{eq:reduction}, and uniqueness follows from the
last-nonzero-coordinate argument. Thus $K_\infty$ is a transversal, and
$K_\infty+\Lambda_\infty$ tiles $\Z^{(\mathbb N)}$.

The inclusion of the first $d$ coordinates induces an injection
$\Z^d/\Lambda_d^\varepsilon\to\Z^{d+1}/\Lambda_{d+1}^\varepsilon$, since a
lattice vector with last coordinate zero uses the last basis vector with
coefficient zero. Hence $\Z^{(\mathbb N)}/\Lambda_\infty$ is an increasing
union of cyclic groups of orders $2^d$, that is, a Pr\"ufer $2$-group. It
is not cyclic, since its elements have unbounded finite orders. For the
binary signs ($\varepsilon_j=-1$, normalized basis $B_d$) the isomorphism
is explicit:
\[
 \Z^{(\mathbb N)}/\Lambda_\infty\ \cong\ \Z[\tfrac12]/\Z,\qquad
 x\longmapsto\sum_{i\ge1}x_i2^{-i}\bmod1,
\]
because $\sum_ix_i2^{-i}=F_d(x)/2^d$ for $x$ supported on the first $d$
coordinates. Under this isomorphism $K_\infty$ corresponds to the finite
binary expansions $0.x_1x_2x_3\ldots$ of the dyadic rationals in $[0,1)$.

Define the infinite-dimensional counter $\ell_\infty(x)$ as the number of
distinct displacements $\rho(x+u)-u$, $u\in K_\infty$. If $x$ is
supported on the first $d$ coordinates, a displacement congruent to $x$
cannot have a nonzero coordinate beyond $d$, since the last nonzero
coordinate of a lattice vector is even. Therefore
\[
 \ell_\infty(x)=\ell_d(x_1,\ldots,x_d)<\infty.
\]
In the binary case, by \eqref{eq:binaryprofile} and \eqref{eq:sternrow},
$\ell_\infty(x)=s(2^d+2^d\theta)$, where $\theta\in[0,1)$ is the residue of
$\sum_ix_i2^{-i}$ modulo $1$ and $d$ is any integer with
$2^d\theta\in\Z$; the right-hand side does not depend on $d$, since
$s(2m)=s(m)$. No counting identity with an infinite sum over $K_\infty$ is
asserted.

%======================================================================
\section*{Declaration of generative AI and AI-assisted technologies}
%======================================================================
During the preparation of this work the authors used ChatGPT Astra Pro~6
(OpenAI) and Claude (Anthropic, model Claude Opus~5.5). ChatGPT Astra
obtained the congruence representation in Theorem~\ref{thm:weights} and,
answering questions posed by the authors, the proofs of the results of
Section~\ref{sec:profiles}, whose statements the authors had anticipated;
the authors suggested a simpler derivation of the weights and the use of the
Fourier transform. Claude merged the authors' drafts into a single
manuscript, edited the text, drew the figures, checked the proofs,
contributed Proposition~\ref{prop:converse}, Remarks~\ref{rem:hnfcount},
\ref{rem:signedweights}, \ref{rem:keplerincw} and~\ref{rem:sternorbit},
Proposition~\ref{prop:rigidity}, Theorem~\ref{thm:boxclass}(b),
Corollary~\ref{cor:boxprofiles}, Remark~\ref{rem:boxorbits} and part of
Appendix~\ref{app:infinite}, and carried out the computations in
Section~\ref{sec:profiles} and Remark~\ref{rem:multisets}. The authors
reviewed and verified all of the content and take full responsibility for
this publication.

%======================================================================

\end{document}